\documentclass[10pt]{amsart}
\usepackage{latexsym,amsmath,amssymb,amscd}
\def\today{\ifcase \month \or
   January \or February \or March \or April \or
   May \or June \or July \or August \or
   September \or October \or November \or December \fi
   \space\number\day , \number\year}

\begin{document}

\makeatletter
\@addtoreset{figure}{section}
\def\thefigure{\thesection.\@arabic\c@figure}
\def\fps@figure{h,t}
\@addtoreset{table}{bsection}

\def\thetable{\thesection.\@arabic\c@table}
\def\fps@table{h, t}
\@addtoreset{equation}{section}
\def\theequation{
\arabic{equation}}
\makeatother

\newcommand{\bfi}{\bfseries\itshape}

\newtheorem{theorem}{Theorem}
\newtheorem{acknowledgment}[theorem]{Acknowledgment}
\newtheorem{claim}[theorem]{Claim}
\newtheorem{conclusion}[theorem]{Conclusion}
\newtheorem{condition}[theorem]{Condition}
\newtheorem{conjecture}[theorem]{Conjecture}
\newtheorem{construction}[theorem]{Construction}
\newtheorem{corollary}[theorem]{Corollary}
\newtheorem{criterion}[theorem]{Criterion}
\newtheorem{definition}[theorem]{Definition}
\newtheorem{example}[theorem]{Example}
\newtheorem{lemma}[theorem]{Lemma}
\newtheorem{notation}[theorem]{Notation}
\newtheorem{problem}[theorem]{Problem}
\newtheorem{proposition}[theorem]{Proposition}
\newtheorem{question}[theorem]{Question}
\newtheorem{remark}[theorem]{Remark}

\numberwithin{theorem}{section}
\numberwithin{equation}{section}

\newcommand{\todo}[1]{\vspace{5 mm}\par \noindent
\framebox{\begin{minipage}[c]{0.95 \textwidth}
#1 \end{minipage}}\vspace{5 mm}\par}

\newcommand{\1}{{\bf 1}}

\newcommand{\hotimes}{\widehat\otimes}

\newcommand{\Alt}{{\rm Alt}\,}
\newcommand{\Ci}{{\mathcal C}^\infty}
\newcommand{\comp}{\circ}
\newcommand{\D}{\text{\bf D}}
\newcommand{\ev}{{\rm ev}}
\newcommand{\id}{{\rm id}}
\newcommand{\ad}{{\rm ad}}
\newcommand{\de}{{\rm d}}
\newcommand{\dist}{{\rm dist}}
\newcommand{\Ad}{{\rm Ad}}
\newcommand{\ie}{{\rm i}}
\newcommand{\End}{{\rm End}\,}
\newcommand{\Gr}{{\rm Gr}_{\rm res}}
\newcommand{\Grr}{{\rm Gr}}
\newcommand{\Hom}{{\rm Hom}\,}
\newcommand{\Ker}{{\rm Ker}\,}
\newcommand{\Lie}{\text{\bf L}}
\newcommand{\lf}{{\rm l}}
\newcommand{\res}{{\rm res}}
\newcommand{\Ran}{{\rm Ran}\,}
\newcommand{\spann}{{\rm span}}
\newcommand{\Tr}{{\rm Tr}\,}

\newcommand{\G}{{\rm G}}
\newcommand{\U}{{\rm U}}
\newcommand{\Z}{{\rm Z}}
\newcommand{\VB}{{\rm VB}}

\newcommand{\Ac}{{\mathcal A}}
\newcommand{\Bc}{{\mathcal B}}
\newcommand{\Cc}{{\mathcal C}}
\newcommand{\Dc}{{\mathcal D}}
\newcommand{\Fc}{{\mathcal F}}
\newcommand{\Hc}{{\mathcal H}}
\newcommand{\Oc}{{\mathcal O}}
\newcommand{\Pc}{{\mathcal P}}
\newcommand{\Qc}{{\mathcal Q}}
\newcommand{\Wc}{{\mathcal W}}
\newcommand{\Xc}{{\mathcal X}}

\newcommand{\Bg}{{\mathfrak B}}
\newcommand{\Jg}{{\mathfrak J}}
\newcommand{\Lg}{{\mathfrak L}}
\newcommand{\Sg}{{\mathfrak S}}
\newcommand{\Xg}{{\mathfrak X}}
\newcommand{\ug}{{\mathfrak u}}
\newcommand{\g}{{\mathfrak g}}
\newcommand{\bg}{{\mathfrak b}}
\newcommand{\hg}{{\mathfrak h}}
\newcommand{\kg}{{\mathfrak k}}

\def\R{{\mathbb{R}}}
\def\C{{\mathbb{C}}}
\def\D{{\mathbb{D}}}
\def\N{{\mathbb{N}}}

\pagestyle{myheadings} \markboth{}{}


\makeatletter
\title{The restricted Grassmannian, Banach Lie-Poisson spaces,
and coadjoint orbits}
\author{Daniel Belti\c t\u a}
\author{Tudor S. Ratiu}
\author{Alice Barbara Tumpach}

\address{Institute of Mathematics ``Simion
Stoilow'' of the Romanian Academy,
P.O. Box 1-764, RO-014700 Bucharest, Romania}
\email{Daniel.Beltita@imar.ro}
\address{Section de Math\'ematiques and Centre  Bernoulli,
\'Ecole Polytechnique F\'ed\'erale de Lausanne,
CH-1015 Lausanne, Switzerland}
\email{Tudor.Ratiu@epfl.ch}
\address{Section de Math\'ematiques,
\'Ecole Polytechnique F\'ed\'erale de Lausanne,
CH-1015 Lausanne, Switzerland}
\email{Alice.Tumpach@epfl.ch}

\begin{abstract}
We investigate some basic questions concerning
the relationship between the restricted Grassmannian and
the theory of Banach Lie-Poisson spaces.
By using universal central extensions of Lie algebras, 
we find that the restricted Grassmannian is symplectomorphic
to symplectic leaves in certain Banach Lie-Poisson spaces,   
and the underlying Banach space can be chosen to be even a Hilbert space.
Smoothness of numerous adjoint and coadjoint orbits
of the restricted unitary group is also established.
Several pathological properties of the restricted algebra
are pointed out.


\end{abstract}

\date{\today}
\makeatother
\maketitle


\section{Introduction}

The present paper is devoted to an investigation of the
relationship between the restricted Grassmannian 
and the recently
initiated theory of Banach Lie-Poisson spaces. 

The restricted Grassmannian 
(whose definition is recalled after Proposition~\ref{leaf} below)
is a quite remarkable infinite-dimensional {\it K\"ahler} manifold 
that plays an important role in many areas of
mathematics and physics. There are many interesting objects related to the restricted Grassmannian, such as: loop groups (see Proposition 8.3.3 in \cite{PS90}), the
coadjoint orbits $\textrm{Diff}^{+}(S^{1})/S^{1}$ and
$\textrm{Diff}^{+}(S^{1})/\textrm{PSU}(1,1)$ of the group of
orientation-preserving diffeomorphisms of the circle (Proposition~6.8.2 in \cite{PS90} and
Proposition~5.3 in \cite{Se81}). 
It is related to the integrable
system defined by the KP hierarchy (see  \cite{SW85}) and to the
fermionic second quantization (see \cite{Wu01}). 
On the other
hand, the notion of a Banach Lie-Poisson space was recently introduced in
\cite{OR03} and is an infinite-dimensional version of the
Lie-Poisson spaces, that is, the Poisson manifolds provided by
dual spaces of finite-dimensional Lie algebras 
(see for instance \cite{OrR04} for the finite-dimensional theory). 
This new class of
infinite-dimensional linear Poisson manifolds is remarkable in
several respects: it includes all the preduals of $W^*$-algebras,
thus establishing a bridge between Poisson geometry and the
theory of operator algebras, and hence it provides links with 
algebraic quantum theories; it interacts in a fruitful way with
the theory of extensions of Lie algebras (see \cite{OR04}); and
finally, there exist large classes of Banach Lie-Poisson spaces
which share with the finite-dimensional Poisson manifolds the
fundamental property that the characteristic distribution is
integrable, the corresponding integral manifolds being in addition
Poisson submanifolds which are symplectic and, in several
important situations, are even {\it K\"ahler} manifolds (see
\cite{BR05}).

We have mentioned here two types of infinite-dimensional
K\"ahler manifolds: the restricted Grassmannians and 
certain symplectic leaves in infinite-dimensional Lie-Poisson
spaces introduced in \cite{OR03}.
This brings us to the first question
addressed in the present paper:

\begin{question}\label{grass_symp_leaf}
Is the restricted Grassmannian a symplectic leaf in a Banach
Lie-Poisson space?
\end{question}

The main result of our paper is essentially affirmative and the precise answer is given in
Section~\ref{MAIN}. Specifically, we shall employ the method of
central extensions to construct a certain Banach Lie-Poisson space
$\widetilde{u}_2$ whose characteristic distribution is integrable
(Theorem~\ref{gr_leaf_in_u2}) and one of the integral
manifolds of this distribution is symplectomorphic to the
connected component $\Gr^0$ of the restricted Grassmannian
(Theorem~\ref{main}). Using a similar method, we realize the
restricted Grassmannian as a symplectic leaf in yet another
Banach Lie-Poisson space, which is the predual to a 1-dimensional
central extension of the restricted Lie algebra $\ug_{\res}$. See
Section~\ref{cen} for a detailed discussion of the Poisson
geometry of this new Banach Lie-Poisson space
$(\tilde{\ug}_{\res})_{*}$.

This second construction is closely related to another area where
the theory of restricted groups interacts with
the theory initiated in \cite{OR03}.
Specifically, we also address the following
question on the predual $(\ug_{\res})_*$ of the restricted Lie algebra:

\begin{question}\label{distrib_integ}
Does the real Banach space $(\ug_{\res})_*$ have a natural structure of
Banach Lie-Poisson space and is its characteristic distribution
integrable?
\end{question}

By the very construction of the Banach Lie-Poisson space
$(\tilde{\ug}_{\res})_{*}$, the predual $(\ug_{\res})_*$ appears
as a Poisson submanifold of $(\tilde{\ug}_{\res})_{*}$ and carries
a natural structure of Banach Lie-Poisson space. Nonetheless, the
answer to the second part of Question~\ref{distrib_integ} turns
out to be much more difficult to give than the one to
Question~\ref{grass_symp_leaf} inasmuch as the restricted algebra
$\Bc_{\res}$ (see Notation~\ref{grass} below) is a dual Banach
$*$-algebra with many pathological properties (summarized
in Section~\ref{pat}): its unitary group is unbounded, its natural
predual is not spanned by its positive cone, and a conjugation
theorem for its maximal Abelian $*$-subalgebras fails to be true.
Despite these unpleasant properties, we show that the characteristic
distribution of $(\ug_{\res})_*$ has numerous smooth integral
manifolds, which are, in particular, {\it smooth} coadjoint
orbits of the restricted unitary group~$\U_{\res}$ (see
Section~\ref{coadj}). For the sake of completeness, a short
section of the paper (Section~\ref{adj}) is devoted to
investigating smoothness of adjoint orbits of $\U_{\res}$.

\begin{notation}\label{grass}
\normalfont We conclude this Introduction by setting up some
notation to be used throughout the paper. In the following, $\Hc$
will denote a separable complex Hilbert space, endowed with a
decomposition $\Hc=\Hc_{+}\oplus\Hc_{-}$ into the orthogonal sum
of two closed infinite-dimensional subspaces. The orthogonal
projection onto $\Hc_{\pm}$ will be denoted by $p_{\pm}$. The
Banach ideal of trace class operators on $\Hc$ will be denoted by
$\Sg_1(\Hc)$ and $\Sg_2(\Hc)$ will denote the Hilbert ideal of
Hilbert-Schmidt operators on $\Hc$. We let $\Bc(\Hc)$ be the
algebra of all bounded linear operators on $\Hc$. We shall also
need the Banach-Lie group of unitary operators on $\Hc$,
$$\U(\Hc)=\{u\in\Bc(\Hc)\mid u^*u=uu^*=\id\}, $$
whose Lie algebra is
$$\ug(\Hc)=\{a\in\Bc(\Hc)\mid a^*=-a\}.$$
Now let us define the following
skew-Hermitian element:
$$
d:=\ie(p_{+}-p_{-})\in\ug(\Hc).
$$
The restricted Banach algebra and the restricted unitary group are
respectively defined as follows:
$$
\begin{array}{l}
\Bc_{\res}=\{a\in\Bc(\Hc)\mid [d,a]\in\Sg_2(\Hc)\}
           =\{a\in\Bc(\Hc)\mid
           \Vert a\Vert_{\res}:=\Vert a\Vert
           +\Vert[d,a]\Vert_2<\infty\},\text{ and}\\
\U_{\res}=\{u\in\U(\Hc)\mid[d,u]\in\Sg_2(\Hc)\}
     =\U(\Hc)\cap\Bc_{\res}.
\end{array}
$$
The Lie algebra of $\U_{\res}$ is the following Banach Lie algebra:
$$
\ug_{\res}=\{a\in\ug(\Hc)\mid[d,a]\in\Sg_2(\Hc)\}
     =\ug(\Hc)\cap\Bc_{\res}.
$$
Let us define the following Banach Lie algebra:
$$
(\ug_{\res})_{*}=\{\rho\in\ug(\Hc)\mid
[d,\rho]\in\Sg_2(\Hc),\;
p_{\pm}\rho|_{\Hc_{\pm}}\in\Sg_1(\Hc_{\pm})\}.
$$
A connected Banach Lie group with Lie algebra $(\ug_{\res})_{*}$
is
$$
\U_{1,2}=\{a\in\U(\Hc)\mid a-\id\in\Sg_2(\Hc),\;
p_{\pm}a|_{\Hc_{\pm}}\in \id+\Sg_1(\Hc_{\pm})\}.
$$
The group $\U_{1}$ and its Lie algebra $\ug_1$ are defined
as follows:
$$
\begin{array}{l}
\U_{1} = \{a\in\U(\Hc)\mid a-\id\in\Sg_1(\Hc)\},\text{ and}\\
\ug_{1}=\ug(\Hc)\cap\Sg_1(\Hc).
\end{array}
$$
Finally, the Hilbert-Lie group $\U_2$ and its Lie
algebra $\ug_2$ are defined by :
$$
\begin{array}{l}
\U_{2} = \{a\in\U(\Hc)\mid a-\id\in\Sg_2(\Hc)\},\text{ and}\\
\ug_{2}=\ug(\Hc)\cap\Sg_2(\Hc).
\end{array}
$$
\end{notation}

\section{The Banach Lie-Poisson space associated to
the universal central extension of $\ug_{\res}$}\label{cen}

In this section we construct a Banach Lie-Poisson space
$(\tilde{\ug}_{\res})_{*}$ whose dual is the universal central
extension of the restricted algebra $\ug_{\res}$.
(See \cite{Ne02b}
for the definition of universal central extension
and Proposition~\ref{cohomology} below for the justification of this fact.)
The Poisson structure of $(\tilde{\ug}_{\res})_{*}$ is defined by
\eqref{bracket_d2} in Proposition \ref{Banach_Poisson}. Let us
first justify the suggestive notation $(\ug_{\res})_{*}$.
\begin{proposition}\label{duality}
The Lie algebra $(\ug_{\res})_{*}$ is a predual of the unitary
restricted algebra $\ug_{\res}$, the duality pairing
$\langle\cdot\,,\cdot\rangle$ being given by
\begin{equation}\label{dualitypairing}
\langle\cdot\,,\cdot\rangle~:(\ug_{\res})_{*}\times\ug_{\res}\to{\mathbb
R},\quad (b,c)\mapsto\Tr(bc).
\end{equation}
\end{proposition}

\begin{proof}
Consider two arbitrary elements
$$a=\begin{pmatrix} a_{++} & a_{+-} \\
                    -a_{+-}^* & a_{--}
    \end{pmatrix}\in\ug_{\res} \quad 
\text{ and } \quad 
\rho=\begin{pmatrix} \rho_{++} & -\rho_{-+}^* \\
                     \rho_{-+} & \rho_{--}
    \end{pmatrix}\in(\ug_{\res})_{*}.$$
Then
\begin{equation}\label{mult}
a\rho=\begin{pmatrix}
       a_{++}\rho_{++}+a_{+-}\rho_{-+} & -a_{++}\rho_{-+}^{*}+a_{+-}\rho_{--} \\
       -a_{+-}^*\rho_{++}+a_{--}\rho_{-+} & a_{+-}^{*}\rho_{-+}^{*}+a_{--}\rho_{--}
    \end{pmatrix},
\end{equation}
hence
\begin{equation}\label{pairing}
\Tr(a\rho)=\Tr(a_{++}\rho_{++})+2\Re\Tr(a_{+-}\rho_{-+})+
\Tr(a_{--}\rho_{--}),
\end{equation}
where $\mathfrak{R}z$ denotes the real part of the complex number $z $. Recall that the bilinear functional
$$
\Bc(\Hc_{\pm})\times\Sg_1(\Hc_{\pm})\to{\mathbb C},\quad
(b,c)\mapsto\Tr(bc),
$$
induces a topological isomorphism of complex Banach spaces
$\left(\Sg_1(\Hc_{\pm})\right)^*\simeq\Bc(\Hc_{\pm})$. It follows
that the trace induces a topological isomorphism of real Banach
spaces
\begin{equation}\label{unitarypairing}
\left(\ug(\Hc_{\pm})\cap\Sg_1(\Hc_{\pm})\right)^{*}\simeq\ug(\Hc_{\pm}).
\end{equation}
Indeed, the $\C$-linearity of the trace implies that for
$b\in\Bc(\Hc_{\pm})$ the following conditions are equivalent:  
$$\bigl(\forall c\in\ug(\Hc_{\pm})\cap\Sg_1(\Hc_{\pm})\bigr)\quad \Tr(bc)=0
\iff
\bigl(\forall c\in\Sg_1(\Hc_{\pm})\bigr)\quad \Tr(bc)=0.
$$
Moreover the condition
$$\bigl(\forall c\in\ug(\Hc_{\pm})\cap\Sg_1(\Hc_{\pm})\bigr)\quad 
\Tr(bc)\in\mathbb{R}
$$
implies
$$\bigl(\forall c\in\ug(\Hc_{\pm})\cap\Sg_1(\Hc_{\pm})\bigr)\quad 
\Tr(b+b^{*})c=0,
$$
hence $b$ belongs to $\ug(\Hc_{\pm})$. 
On the other hand, the
duality pairing of complex Hilbert spaces
$$
\Sg_2(\Hc_{-},\Hc_{+})\times\Sg_2(\Hc_{+},\Hc_{-})\to{\mathbb
\C},\quad (b,c)\mapsto\Tr(bc),
$$
induces a duality pairing of the underlying real Hilbert spaces by
\begin{equation}\label{s2r}
\Sg_2(\Hc_{-},\Hc_{+})\times\Sg_2(\Hc_{+},\Hc_{-})\to{\mathbb
\R},\quad (b,c)\mapsto\Re\,\Tr(bc).
\end{equation}
In view of formula \eqref{pairing}, we conclude that the trace
induces a topological isomorphism of real Banach spaces
$$((\ug_{\res})_*)^*\simeq\ug_{\res}.$$
That is, $(\ug_{\res})_*$ is indeed a predual to $\ug_{\res}$, the
duality pairing being induced by
\eqref{unitarypairing} and \eqref{s2r}.
\end{proof}

\begin{definition}
\normalfont 
We define the Banach Lie algebra $\tilde{\ug}_{\res}$
as the central extension of $\ug_{\res}$ with continuous
two-cocycle $s$  given by
\begin{equation}\label{schwinger}
s(A, B) := \Tr (A[d, B]),
\end{equation}
for all $A,B \in \ug_{\res}$. That is, $\tilde{\ug}_{\res}$ is
the Banach algebra $\ug_{\res} \oplus \R$ endowed with the bracket
$[\cdot, \cdot]_{d}$ defined by
\begin{equation}\label{bracket u^d}
\left[(A, a),(B, b)\right]_{d} = \left([A, B], -s(A, B)\right).
\end{equation}
\end{definition}

\begin{remark}\label{non-triv}
\normalfont 
Note that by the very definition of $\ug_{\res}$, one
has $[d, \ug_{\res}] \subset (\ug_{\res})_{*}$. It follows from
the duality pairing \eqref{dualitypairing}, that $s$ is
well-defined by \eqref{schwinger}. 
To see that $s$ defines a
two-cocycle on $\ug_{\res}$, let us remark that $s$ is
($2\ie$)-times the Schwinger term of \cite{Wu01}. It follows from
Corollary~II.12 in the aforementioned work that $s$ defines a
non-trivial element in the second continuous Lie algebra
cohomology space $H^2(\ug_{\res}, \mathbb{R})$. 
The corresponding
$\U(1)$-extension of the unitary restricted group $\U_{\res}$ is
isomorphic to the $\U(1)$-extensions ${\U}_{\res}^{\sim}$ and
$\widehat{\U}_{\res}$ of $\U_{\res}$ constructed in \cite{Wu01}.
\end{remark}

\begin{proposition}\label{cohomology}
The cohomology class $[s]$ is a generator of
the continuous Lie algebra cohomology space~$H^2(\ug_{\res}, \mathbb{R})$.
\end{proposition}

\begin{proof}
According to Proposition~I.11 in \cite{Ne02a}, the second
continuous Lie algebra cohomology space $H^2(\Bc_{\res},
\mathbb{C})$ of the restricted Lie algebra $\Bc_{\res}$ is 1-dimensional.
Note that a continuous $\R$-valued $2$-cocycle $v$ on $\ug_{\res}$
extends by $\C$-linearity to a continuous $\C$-valued $2$-cocycle
$v^{\C}$ on the complex Lie algebra $\Bc_{\res}$. 
The cocycle
$v^{\C}$ is a coboundary if and only if there exists a continuous
linear map  $\alpha\colon\Bc_{\res}\to\C$ such that
$v^{\C}(x, y) = \alpha\left([x,y]\right)$ for every $x,y
\in \Bc_{\res}$. But since $v^{\C}$ restricts to the $\R$-valued
$2$-cocycle $v$ on $\ug_{\res}$, this is the case if and only if
there exists $\beta :=\Re\alpha\colon\ug_{\res}\to\R$ such
that $v(x, y) = \beta\left([x,y]\right)$ for every $x,y
\in \ug_{\res}$. It follows that the extension $v^{\C}$ is a
coboundary on $\Bc_{\res}$ if and only if $v$ is a coboundary on
$\ug_{\res}$. 
Consequently, there is a natural linear injection of
$H^{2}(\ug_{\res}, \R)$ into $H^{2}(\Bc_{\res}, \C)$. 
Since $s$
defines a non-trivial element in $H^{2}(\ug_{\res}, \R)$ (see
Remark~\ref{non-triv}) and $\textrm{dim}_{\C}H^2(\Bc_{\res},
\mathbb{C})= 1$, it follows that
$\textrm{dim}_{\R}H^{2}(\ug_{\res}, \R)=1$ and thus $H^{2}(\ug_{\res},
\R)$ is generated by $s$.
\end{proof}

\begin{proposition}\label{Banach_Poisson}
The Banach space $(\tilde{\ug}_{\res})_{*}$ is a Banach
Lie-Poisson space for the Poisson bracket
\begin{equation}\label{bracket_d2}
\{f,g\}_{d}(\mu, \gamma) := \langle \mu,  \left[ D_{\mu}f(\mu),
D_{\mu}g(\mu) \right] \rangle - \gamma s(D_{\mu}f, D_{\mu}g)
\end{equation}
where $f, g \in \Ci((\tilde{\ug}_{\res})_{*})$, $(\mu,
\gamma)$ is an arbitrary element in $(\tilde{\ug}_{\res})_{*}$,
and $D_{\mu}$ denotes the partial Fr\'echet derivative with
respect to $\mu \in (\ug_{\res})_{*}$.
\end{proposition}

The pairing in equation \eqref{bracket_d2} is the duality pairing
defined by \eqref{dualitypairing}. We will denote by
$\langle\cdot\,, \cdot\rangle_{d}$ the duality pairing between
$(\tilde{\ug}_{\res})_{*} = (\ug_{\res})_{*} \oplus \R$ and
$\tilde{\ug}_{\res} = \ug_{\res} \oplus \R$ given by
$$
\left\langle (\mu, \gamma), (A, a) \right\rangle_{d} = \langle
\mu, A\rangle + \gamma a.
$$

\begin{proof}[Proof of Proposition~\ref{Banach_Poisson}]
By Theorem~4.2 in \cite{OR03}, the Banach space
$(\tilde{\ug}_{\res})_{*}$ is a Banach Lie-Poisson space if and
only if its dual $\tilde{\ug}_{\res}$ is a Banach Lie algebra
satisfying $\ad_{x}^{*}(\tilde{\ug}_{\res})_{*} \subset
(\tilde{\ug}_{\res})_{*} \subset (\tilde{\ug}_{\res})^{*}$ for all
$x \in \tilde{\ug}_{\res}$. 
The fact that $\tilde{\ug}_{\res}$ is
a Banach Lie algebra follows directly from the continuity of $s$
and from the $2$-cocycle identity which implies the Jacobi
identity of $[\cdot,\cdot]_{d}$. 
To see that the coadjoint action
of $\tilde{\ug}_{\res}$ preserves the predual
$(\tilde{\ug}_{\res})_{*}$, note that for every $(A, a),(B, b)
\in \tilde{\ug}_{\res}$ and every $(\mu, \gamma) \in
(\tilde{\ug}_{\res})_{*}$, one has
$$
\begin{aligned}
\langle -\ad_{(A, a)}^{*}(\mu, \gamma), (B, b)\rangle_{d} := &
  \langle (\mu, \gamma), \ad_{(A, a)}(B, b)\rangle_{d}
 =  \langle (\mu,
\gamma), [(A, a), (B, b)]_{d}\rangle_{d}\\
= & \langle (\mu, \gamma), \left([A, B], - s(A, B)\right)
\rangle_{d}
 = \Tr \mu [A, B] - \gamma \Tr A[d, B]\\
= & \Tr \mu [A, B]+ \gamma \Tr [d, A] B   =  \langle
(-\ad^{*}(A)(\mu) + \gamma [d, A], 0), (B, b) \rangle_{d}.
\end{aligned}
$$
Since
\begin{equation}\label{ideals}
[(\ug_{\res})_*,\ug_{\res}]\subseteq(\ug_{\res})_*,
\end{equation}
and
\begin{equation}\label{derivation}
[d, \ug_{\res}] \subset (\ug_{\res})_{*},
\end{equation}
we conclude that $-\ad^{*}(A)(\mu) + \gamma [d, A]$ belongs to
$(\ug_{\res})_{*}$ for every $A \in \ug_{\res}$. Hence the predual
$(\tilde{\ug}_{\res})_{*}$ is preserved by the coadjoint action.
Referring again to Theorem~4.2 in \cite{OR03}, it follows that the
Poisson bracket of $f$, $g \in \Ci((\tilde{\ug}_{\res})_{*})$ is
given by
$$
\{f, g\}_{d}(\mu, \gamma) = \langle (\mu, \gamma),
[Df(\mu,\gamma), Dg(\mu, \gamma)] \rangle_{d}.
$$
Denoting respectively by $D_{\mu}$ and $D_{\gamma}$ the partial
Fr\'echet derivatives with respect to $\mu \in (\ug_{\res})_{*}$
and $\gamma \in \R$, one has
$$\begin{aligned}
\{f, g\}_{d}(\mu, \gamma) &= \langle(\mu, \gamma),
\left[\left(D_{\mu}f, D_{\gamma}f\right), \left(D_{\mu}g,
D_{\gamma}g\right)\right]_{d}
\rangle_{d}\\
&= \langle (\mu, \gamma), \left([D_{\mu}f, D_{\mu}g],
-s(D_{\mu}f,
D_{\mu}g)\right)\rangle_{d}\\
&= \langle \mu, [D_{\mu}f, D_{\mu}g]\rangle - \gamma s(D_{\mu}f,
D_{\mu}g),
\end{aligned}
$$
and this ends the proof.
\end{proof}

\begin{remark}\label{Xh}
\normalfont
By Theorem~4.2 in \cite{OR03}, it follows that the
Hamiltonian vector field associated to a smooth function $h$ on
$(\ug_{\res})_{*}$ is given by
\begin{equation}\label{hamiltonian}
X_{h}(\mu, \gamma) = -\textrm{ad}^{*}_{(D_{\mu}h,
D_{\gamma}h)}(\mu,
\gamma)
= \bigl(-\ad^{*}_{D_{\mu}h}\mu - \gamma [D_{\mu}h, d], 0\bigr).
\end{equation}
\end{remark}

\begin{remark}
\normalfont 
Note that, for each $\gamma \in \R$, $(\ug_{\res})_{*}
\oplus \{\gamma\} $ is a Poisson submanifold of
$(\tilde{\ug}_{\res})_{*}$  for the following  Poisson bracket on
the first factor
$$
\{f,g\}_{d, \gamma}(\mu) := \langle \mu, \left[ D_{\mu}f(\mu),
D_{\mu}g(\mu) \right] \rangle - \gamma s(D_{\mu}f, D_{\mu}g).
$$
\end{remark}

\begin{remark}
\normalfont The central extension $(\tilde{\ug}_{\res})_{*}$ of
the Banach Lie-Poisson space $(\ug_{\res})_{*}$ is a particular
example of the extensions of Banach Lie-Poisson spaces constructed
in \cite{OR04}. Indeed formula \eqref{bracket_d2} for the bracket
of two functions on $(\tilde{\ug}_{\res})_{*}$ can be
alternatively deduced from the general formula~(5.6) in
Theorem~5.2 of \cite{OR04}, with $\mathfrak{c}=\mathbb{R}$,
$\mathfrak{a}=(\ug_{\res})_*$, $\varphi=0$ and $\omega=-s$. The
pairing in the second term of the right hand side of (5.6),
Theorem~5.2, \cite{OR04}, is, in this special case, just the
pairing between the real line and its dual given by multiplication
of real numbers (the element $c\in\mathfrak{c}$ is $\gamma$), and
the bracket of partial derivatives of the functions $f$ and $g$
with respect to $c$ vanishes since $\mathbb{R}$ is commutative.
\end{remark}

\begin{proposition}\label{affine_coadjoint}
The  unitary group $\U_{\res}$ acts on the Poisson manifold
$(\ug_{\res})_{*} \oplus \{\gamma\} \subset
(\tilde{\ug}_{\res})_{*}$ by affine coadjoint action as follows.
For  $g \in \U_{\res}$,
$$
\begin{array}{lcl}
g\cdot(\mu, \gamma) := \left(\Ad^{*}(g^{-1})(\mu) -
\gamma\sigma(g), \gamma\right)
\end{array}
$$
where $\mu \in (\ug_{\res})_{*}$, $\gamma\in\mathbb{R}$, and where
$$
\begin{aligned}
\sigma\colon \U_{\res} & \to (\ug_{\res})_{*},\\
                     g & \mapsto  gdg^{-1} - d.
\end{aligned}
$$
\end{proposition}

\begin{proof}
Let us verify that for every $g \in \U_{\res}$ we have
$g\,d\,g^{-1} -d\in(\ug_{\res})_{*}$.
Consider the block
decomposition of $g$ with respect to the direct sum
$\Hc=\Hc_{+}\oplus\Hc_{-}$
$$
g =\begin{pmatrix} g_{++} & g_{+-} \\
                    g_{-+} & g_{--}
    \end{pmatrix}\in\U_{\res}.
$$
 One has
\begin{equation}\label{product}
\begin{pmatrix} g_{++} & g_{+-} \\
                    g_{-+} & g_{--}
    \end{pmatrix}
\begin{pmatrix} \ie & 0 \\
                    0 & -\ie
    \end{pmatrix}
    \begin{pmatrix} g_{++}^{*} & g_{-+}^{*} \\
                    g_{+-}^{*} & g_{--}^{*}
    \end{pmatrix}
    =
\begin{pmatrix}
\ie g_{++}g_{++}^{*}-\ie g_{+-}g_{+-}^{*}
      & \ie g_{++}g_{-+}^{*}-\ie g_{+-}g_{--}^{*} \\
\ie g_{-+}g_{++}^{*}-\ie g_{--}g_{+-}^{*}
      & \ie g_{-+}g_{-+}^{*}-\ie g_{--}g_{--}^{*}.
    \end{pmatrix}
\end{equation}
Since $g_{\pm\mp}$ belongs to $\Sg_{2}(\Hc_{\mp},\Hc_{\pm})$, the
off-diagonal blocks of the right hand side are in
$\Sg_{2}(\Hc_{\pm},\Hc_{\mp})$. Further, since
$$
\begin{pmatrix} g_{++} & g_{+-} \\
                    g_{-+} & g_{--}
    \end{pmatrix}
    \begin{pmatrix} g_{++}^{*} & g_{-+}^{*} \\
                    g_{+-}^{*} & g_{--}^{*}
    \end{pmatrix}
    =
\begin{pmatrix}
g_{++}g_{++}^{*} + g_{+-}g_{+-}^{*} & g_{++}g_{-+}^{*} + g_{+-}g_{--}^{*} \\
g_{-+}g_{++}^{*} + g_{--}g_{+-}^{*} & g_{-+}g_{-+}^{*} + g_{--}g_{--}^{*}
    \end{pmatrix}
=
\begin{pmatrix} \id & 0 \\
                  0 & \id
    \end{pmatrix},
$$
and since $\Sg_{2}\cdot\Sg_{2} \subset \Sg_{1}$,  one has
$$g_{++}g_{++}^{*} = \id-g_{+-}g_{+-}^{*} \in \id+\Sg_{1}(\Hc_{+})$$
and
$$g_{--}g_{--}^{*} = \id-g_{-+}g_{-+}^{*} \in \id+\Sg_{1}(\Hc_{-}).$$
Consequently,
$$
g_{++}g_{++}^{*}-g_{+-}g_{+-}^{*}\in \id+\Sg_{1}(\Hc_{+})
$$
and
$$ g_{-+}g_{-+}^{*} -g_{--}g_{--}^{*}\in
-\id+\Sg_{1}(\Hc_{-}).
$$
Moreover, it is clear that the result of the multiplication
(\ref{product}) is skew-symmetric.
Hence for all $g \in
\U_{\res}$ we have $g\,d\,g^{-1} - d\in (\ug_{\res})_{*}$.

Denoting by $\textrm{Aff}\left((\ug_{\res})_{*}\right)$ the affine
group of transformations of $(\ug_{\res})_{*}$, it remains to show
that
$$
\begin{aligned}
(\Ad^*,\gamma\sigma)\colon \U_{\res} &\to
\textrm{Aff}((\ug_{\res})_{*}) = \textrm{GL}((\ug_{\res})_{*})
\rtimes (\ug_{\res})_{*}\\
g &\mapsto(\Ad^{*}(g^{-1}), \gamma\sigma(g))
\end{aligned}
$$
is a group homomorphism. For this, we have to check that $
\gamma\sigma(g_{1}g_{2}) = \Ad^{*}(g_{1}^{-1})\gamma\sigma(g_{2})
+ \gamma\sigma(g_{1})$ for all $g_{1}$, $g_{2}$ in $\U_{1,2}$ (see
\cite{Ne00}).
In fact
$$
\begin{aligned}
\sigma(g_{1}g_{2})
&= g_{1}g_{2}\,d\,g_{2}^{-1}g_{1}^{-1} - d =
  g_{1}\left(g_{2}\,d\,g_{2}^{-1} - d\right)g_{1}^{-1}
   +(g_{1}\,d\,g_{1}^{-1}- d) \\
&=\Ad^{*}(g_{1}^{-1})\left(\sigma(g_{2})\right) + \sigma(g_{1}),
\end{aligned}
$$
and this ends the proof.
\end{proof}

\begin{proposition}\label{isotropy}
The isotropy group of $(0,
\gamma)\in(\ug_{\res})_{*}\oplus\{\gamma\}$ for the
$\U_{\res}$-affine coadjoint action  is a Lie subgroup of
$\U_{\res}$.
\end{proposition}

\begin{proof}
An element $X$ in the Lie algebra $\ug_{\res}$ of $\U_{\res}$
induces by infinitesimal affine coadjoint action on
$(\ug_{\res})_{*}\oplus\{\gamma\}$ the following vector field
$$
\begin{aligned}
X\cdot(\mu, \gamma) :=&
\frac{d}{dt}\left[\exp(tX)\cdot(\mu,\gamma)\right]_{t=0}\\
=& \Bigl(\frac{d}{dt}\left[\Ad^{*}(\exp(-tX))(\mu) -
\gamma\sigma(\exp(tX))\right]_{t=0}, 0\Bigr)\\
=&\left(-\ad^{*}_{X}(\mu) - \gamma[X, d], 0\right).
\end{aligned}
$$
By definition, the Lie algebra of the isotropy group of $(\mu,
\gamma)$ is
$$
\ug_{(\mu, \gamma)} := \left\{ X \in \ug_{\res}\mid
-\ad^{*}(X)(\mu) + \gamma[X, d]=0 \right\}
$$
The proposition is trivial when $\mu$ and $\gamma$ vanish. For
$\mu = 0$ and $\gamma\neq0$, the Lie algebra $\ug_{(0, \gamma)}$
consist of all elements of $\ug_{\res}$ which commute with $d$.
Hence, for $\gamma\neq0$, $\ug_{(0, \gamma)}
=\ug(\Hc_{+})\oplus\ug(\Hc_{-})$. A topological complement to
$\ug_{(0, \gamma)}$ in $\ug_{\res}$ is
$\mathfrak{m}:=\ug(\Hc)\cap\left(\Sg_{2}(\Hc_{+},
\Hc_{-})\oplus\Sg_{2}(\Hc_{-},\Hc_{+})\right)$.
\end{proof}

\begin{proposition}\label{leaf}
The smooth affine coadjoint orbits of $\U_{\res}$ are tangent to the
characteristic distribution of the Poisson manifold
$(\tilde{\ug}_{\res})_{*}$.
\end{proposition}

\begin{proof}
 By the proof of Proposition~\ref{isotropy}, the image of the
 differential of the orbit map is
$$
 \ug_{\res}\cdot(\mu, \gamma) =
\left\{ \bigl(-\ad^{*}_{X}(\mu) - \gamma[X, d],0\bigr)\mid X \in
\ug_{\res}\right\}.
$$
By Remark~\ref{Xh}, the characteristic space at $(\mu,
\gamma)\in(\tilde{\ug}_{\res})_{*}$ is
$$
P(\mu, \gamma) 
= \{ X_{h}(\mu) =
\bigl(-\ad^{*}_{D_{\mu}h}\mu -\gamma[D_{\mu}h, d], 0\bigr)\mid 
h\in\Ci((\ug_{\res})_{*})\}
= \left\{\left(-\ad^{*}_{X}\mu - \gamma[X, d], 0\right)\mid
X\in\ug_{\res}\right\}.
$$
Thus the assertion follows.
\end{proof}

The {\it restricted Grassmannian} $\Gr$ is defined as the set of
subspaces $W$ of the Hilbert space $\Hc$ such that the orthogonal
projection from $W$ to $\Hc_{+}$ (respectively to $\Hc_{-}$) is a
Fredholm operator (respectively a Hilbert-Schmidt operator).
It follows from Propositions 7.1.2~and~7.1.3 in \cite{PS90} that $\Gr$ is a
Hilbert manifold and a homogeneous space under the natural
action of $\U_{\res}$.
According to Proposition~II.2 in \cite{Wu01},
the connected components of $\U_{\res}$ are the sets
$$
\U_{\res}^{k} = \left\{ \begin{pmatrix} U_{++} & U_{+-} \\
                    U_{-+} & U_{--}
    \end{pmatrix}\in \U_{\res}~~|~~\textrm{index}(U_{++})=k \right\}
\quad \text{for} \quad k\in{\mathbb Z}.
$$
The pairwise disjoint sets
$$
\Gr^k = \left\{ W \in \Gr~~|~~\textrm{index}(p_+|_W\colon W
\rightarrow H_{+})=k \right\},
\quad k\in{\mathbb Z}
$$
are the images of the connected components of $\U_{\res}$ by the
continuous projection $\U_{\res} \rightarrow \Gr =
\U_{\res}/\left(\U(\Hc_{+})\times\U(\Hc_{-})\right)$, and thus
they are the connected components of $\Gr$. In particular, the
connected component of $\Gr$ containing $\Hc_+$ is $\Gr^0$. The
K\"ahler structure of the restricted Grassmannian is defined in
\cite{PS90}, Section 7.8. According to the convention in
\cite{PS90}, the K\"ahler form $\omega_{\textrm{Gr}}$ of $\Gr$ is
the $\U_{\res}$-invariant $2$-form whose value at $\Hc_{+}$ is
given by
\begin{equation}\label{omegagr}
\omega_{\textrm{Gr}}(X, Y) = 2\Im\Tr(X^{*}Y),
\end{equation}
where $X, Y \in \Sg_{2}(\Hc_{+}, \Hc_{-})\simeq
T_{\Hc_{+}}\Gr$ and $\mathfrak{I}z$ denotes the imaginary part of $z \in \mathbb{C}$.. Equivalently, $\omega_{\textrm{Gr}} $ is the quotient of
the following
real-valued anti-symmetric bilinear form $\Omega_{\textrm{Gr}}$ on
$\ug_{\res}$ which vanishes on $\ug(\Hc_+)\oplus\ug(\Hc_-)$ and is
invariant under the $\U(\Hc_+)\times\U(\Hc_-)$-action (see Corollary III.8 in
\cite{Wu01})~:
\begin{equation}\label{omegagrhom}
\Omega_{\textrm{Gr}}(A, B) = -2 s(A, B)
\end{equation}
where $A$ and $B$ belongs to $\ug_{\res}$.
In this correspondence,
an element $A =\begin{pmatrix} A_{++} & A_{-+}^* \\
                    A_{-+} & A_{--}
    \end{pmatrix}$ in $\ug_{\res}$ is identified with the vector
    $X = A_{-+}$ in $\Sg_{2}(\Hc_{+}, \Hc_{-})\simeq
T_{\Hc_{+}}(\Gr)$.

\begin{proposition}\label{leaf=strong}
For every $\gamma\neq0$, the connected components of the
$\U_{\res}$-affine coadjoint orbit $\mathcal{O}_{(0,\gamma)}$ of
$(0, \gamma)\in(\ug_{\res})_{*}\oplus\{\gamma\}$ are
 strong symplectic leaves in the Banach Lie-Poisson space
$(\tilde{\ug}_{\res})_{*}$.
\end{proposition}

\begin{proof}
This follows from the general results given in Theorems~7.3, 7.4
and 7.5 in \cite{OR03}, from the above Propositions
\ref{Banach_Poisson}, \ref{affine_coadjoint}, \ref{isotropy},
\ref{leaf}, and from the fact that the characteristic subspace
$P(0, \gamma)$ is the space
$\mathfrak{m}:=\ug(\Hc)\cap\left(\Sg_{2}(\Hc_{+},
\Hc_{-})\oplus\Sg_{2}(\Hc_{-},\Hc_{+})\right)$, which is closed in
$(\tilde{\ug}_{\res})_{*}$ (we identify the subspace
$\mathfrak{m}$ in $(\ug_{\res})_{*}$ with the subspace
$\mathfrak{m}\oplus\{0\}$ in $(\tilde{\ug}_{\res})_{*}$).
\end{proof}

\begin{theorem}\label{gr_leaf_in_predual}
The connected components of the restricted Grassmannian are strong
symplectic leaves in the Banach Lie-Poisson space
$(\tilde{\ug}_{\res})_{*}$. More precisely, for every
$\gamma\neq0$, the $\U_{\res}$-affine coadjoint orbit
$\mathcal{O}_{(0,\gamma)}$ of 
$(0,\gamma)\in(\ug_{\res})_{*}\oplus\{\gamma\}$ is isomorphic to the
restricted Grassmannian $\Gr$ via the application
$$
\begin{aligned}
\Phi_{\gamma}\colon\Gr & \to\mathcal{O}_{(0,\gamma)}\\
                     W & \mapsto 2\ie\gamma(p_{W}-p_{+}),
\end{aligned}
$$
where $p_{W}$ denotes the orthogonal projection on $W$. The
pull-back by $\Phi_{\gamma}$ of the symplectic form on
$\mathcal{O}_{(0,\gamma)}$ is $(\gamma/2)$-times the
symplectic form $\omega_{\Grr}$ on $\Gr$.
\end{theorem}

\begin{proof}
An element $\rho$ of the affine coadjoint orbit
$\mathcal{O}_{(0,\gamma)}$ of $~(0, \gamma)$ is of the form
$$
\rho = \gamma(g\,d\,g^{-1} - d) =
2\ie\gamma(g\,p_{+}\,g^{-1}-p_{+}),
$$
for some $g \in \U_{\res}$. By Corollary III.4 ii) in \cite{Wu01},
$\Phi_{\gamma}$ is bijective for $\gamma\neq0$. Since the manifold
structure of the orbit $\mathcal{O}_{(0,\gamma)}$ is induced by
the identification $\mathcal{O}_{(0,\gamma)}=
\U_{\res}/\left(\U(\Hc_{+})\times\U(\Hc_{-})\right)$, it follows
from Corollary III.4 i) in \cite{Wu01} that $\Phi_{\gamma}$ is a
diffeomorphism. The symplectic form $\omega_{\mathcal{O}}$ on
$\mathcal{O}_{(0,\gamma)}$ is the $\U_{\res}$-invariant symplectic
form whose value at $(0,\gamma)\in \mathcal{O}_{(0,\gamma)}$ is
the given by
$$
\omega_{\mathcal{O}}\left(0,\gamma\right)\left(X_{f}(0,\gamma),
X_{g}(0,\gamma)\right) = \{f\,,g\}_{d}(0\,,\gamma),
$$
where $f$ and $g$ are any smooth function on $(\ug_{\res})_{*}$.
 Using formula
\eqref{hamiltonian} and \eqref{bracket_d2}, it then follows that
$$
\omega_{\mathcal{O}}\left(0,\gamma\right)\left(\gamma[D_{\mu}f,
d]\,, \gamma[D_{\mu}g, d]\right) = -\gamma s(D_{\mu}f\,,D_{\mu}g).
$$
Hence for every $A, B \in \ug_{\res}$, one has~:
$$
\omega_{\mathcal{O}}\left(0,\gamma\right)\left(\gamma[A, d]\,,
\gamma[B, d]\right) = -\gamma s(A\,,B) =
\frac{\gamma}{2}\Omega_{\textrm{Gr}}(A, B).
$$
It follows that the real-valued anti-symmetric bilinear form on
$\ug_{\res}$ corresponding to the symplectic form
$\omega_{\mathcal{O}}$ on $\mathcal{O}_{(0,\gamma)} =
\U_{\res}/\left(\U(\Hc_+)\times\U(\Hc_-)\right)$ equals
$\frac{\gamma}{2}\Omega_{\textrm{Gr}}$ (where the latter
identification is given by the orbit map), and this ends the
proof.
\end{proof}

\begin{remark}\label{extensions}
\normalfont We refer to the paper \cite{OR04} for additional
information on the relationship between the Banach Lie-Poisson
spaces and the theory of Lie algebra extensions.
\end{remark}

\section{Coadjoint orbits of the restricted unitary group}\label{coadj}

This section includes some partial answers to
Question~\ref{distrib_integ}.
The main difficulty is to show
that the isotropy group of an element in the predual
$(\ug_{\res})_*$ is a Lie subgroup of $\U_{\res}$, or equivalently
that its Lie algebra is complemented in $\ug_{\res}$. Using the
averaging method developed in \cite{Ba90} and \cite{BP05}
for constructing closed
complements, we will be able to show that the $\U_{\res}$-coadjoint
orbit of every element $\rho \in(\ug_{\res})_*$ which commutes with
$d$ is a smooth manifold and that its connected components are
symplectic leaves of the characteristic distribution (see
Proposition \ref{block_diagonal}). 
It follows that the same
conclusion holds for every element $\rho\in(\ug_{\res})_*$ which
is $\U_{\res}$-conjugate to an element commuting to $d$, or
equivalently to a diagonal operator with respect to a Hilbert
basis compatible with the eigenspaces of $d$.
Nevertheless, the set of
elements with the latter property is far from being equal to the whole
$(\ug_{\res})_{*}$.
Recall that in finite dimensions, every
element in the Lie algebra $\ug(n)$ of the unitary group $\U(n)$
is $\U(n)$-conjugate to a diagonal matrix, or, in
other words, $\U(n)$ acts transitively on the set of Cartan
subalgebras of $\ug(n)$. 
This is no longer true in the
infinite-dimensional case (see subsection~\ref{cartan}). 
It is a
difficult question to decide whether a given operator $\rho$ in
$(\ug_{\res})_*$ or $\ug_{\res}$ has the good property of being
$\U_{\res}$-conjugate to a diagonal operator.
In Propositions \ref{finite_rank} and \ref{hinkkanen},
we give some concrete
criteria to check that property.

\begin{conjecture}\label{smoothness}
The real Banach space $(\ug_{\res})_*$ has a natural structure
of Banach Lie-Poisson space and its characteristic distribution
is integrable.
\end{conjecture}

\begin{remark}\label{attempt}
\normalfont
It is clear that
$$(\ug_{\res})_*\hookrightarrow\ug_{\res}$$
with a continuous inclusion map.
On the other hand,
it follows at once by the multiplication formula~\eqref{mult}
that
\begin{equation}\label{ideal}
[(\ug_{\res})_*,\ug_{\res}]\subseteq(\ug_{\res})_*,
\end{equation}
which implies that the predual $(\ug_{\res})_*$ is left invariant
by the coadjoint representation of the Banach Lie algebra
$\ug_{\res}$. Now the results of \cite{OR03} imply the following
two facts:
\begin{itemize}
\item[$\bullet$] The predual Banach space $(\ug_{\res})_*$
has a natural structure of Banach Lie-Poisson space.
\item[$\bullet$] If $\rho\in(\ug_{\res})_*$ has the property that
the corresponding isotropy group
$$
\U_{\res,\rho}:=\{u\in\U_{\res}\mid u\rho u^{-1}=\rho\}
$$
is a Banach Lie subgroup of $\U_{\res}$, then the coadjoint orbit
$\Oc_\rho$ is an integral manifold of the characteristic
distribution of $(\ug_{\res})_*$. Moreover, $\Oc_\rho$ is a weakly
symplectic manifold when equipped with the orbit symplectic
structure.
\end{itemize}
Thus, the desired conclusion will follow as soon as we prove that
the isotropy group $\U_{\res,\rho}$ of any $\rho\in(\ug_{\res})_*$
is a Banach Lie subgroup of $\U_{\res}$.

The Lie algebra of $\U_{\res,\rho}$ is given by
$$\ug_{\res,\rho}=\{a\in\ug_{\res}\mid a\rho=\rho a\}
=\{a\in\ug_{\res}\mid (\forall t\in{\mathbb R})\quad \alpha_t(a)=a\},$$
 where
$$\alpha\colon{\mathbb R}\to \Bc(\ug_{\res}),\quad
\alpha(t)b:=\alpha_t(b):=\exp(t\rho)\cdot b\cdot\exp(-t\rho).$$ It
is clear that $\alpha$ is a group homomorphism. Moreover, since
$\rho\in(\ug_{\res})_*\subseteq\ug_{\res}$ and the adjoint action
of the Banach Lie group $\U_{\res}$ is continuous, it follows that
$\alpha\colon{\mathbb R}\to \Bc(\ug_{\res})$ is norm continuous.

On the other hand, it follows by~\eqref{ideal} that
\begin{equation}\label{weak}
(\forall t\in{\mathbb R})\quad\alpha_t((\ug_{\res})_*)\subseteq(\ug_{\res})_*,
\end{equation}
since $\rho\in(\ug_{\res})_*$.
Then the concrete form of the duality pairing between
$(\ug_{\res})_*$ and $\ug_{\res}$ (see~\eqref{pairing})
shows that
\begin{equation}\label{duals}
(\forall t\in{\mathbb R})\quad
(\alpha_t|_{(\ug_{\res})_*})^*=\alpha_{-t},
\end{equation}
and in particular each operator
$\alpha_t\colon\ug_{\res}\to\ug_{\res}$ is weak$^*$ continuous.

Now a complement to $\ug_{\res,\rho}$ in $\ug_{\res}$ can be constructed
by the averaging technique over the amenable group
$({\mathbb R},+)$ provided one  has
$\sup\limits_{t\in{\mathbb R}}\Vert\alpha_t\Vert<\infty$.
(Some references for the aforementioned averaging technique are
\cite{Ba90}, the proof of Proposition~3.4 in \cite{BR05},
and \cite{BP05}.)

Additionally we note that since for every operator $T:X\to Y$
between the Banach spaces $X$ and $Y$ the norm of $T$ equals the
norm of its dual $T^*$,  it is enough to estimate uniformly the
norm of $\alpha_t$ restricted to the predual $({\mathfrak
u}_{\res})_{*}$. This restriction is an adjoint action of the
group corresponding to the predual.
\end{remark}

\begin{proposition}\label{block_diagonal}
If $\rho\in(\ug_{\res})_*$ and $[d,\rho]=0$, then the coadjoint
isotropy group of $\rho$ is a Banach Lie subgroup of $\U_{\res}$
and the connected components of the corresponding
$\U_{\res}$-coadjoint orbit $\Oc_\rho$ are smooth leaves of the
characteristic distribution of $(\ug_{\res})_*$.
\end{proposition}

\begin{proof}
According to Remark~\ref{attempt} it suffices to show that
$\sup\limits_{t\in{\mathbb R}}\Vert\alpha_t\Vert<\infty$.
The hypothesis  $[d,\rho]=0$ shows that
$\rho$ preserves $\mathcal{H}_{+}$ and $\mathcal{H}_{-}$, that is
$$
\rho=\begin{pmatrix} \rho_{++} & 0 \\
                     0 & \rho_{--}
    \end{pmatrix}\in(\ug_{\res})_{*}.
$$
An element $b \in ({\mathfrak u}_{\res})_{*}$ with block
decomposition with respect to the direct sum
$\Hc=\Hc_{+}\oplus\Hc_{-}$
$$
b=\begin{pmatrix} b_{++} & b_{+-} \\
                    b_{-+} & b_{--}
    \end{pmatrix}
$$
is the sum of an element
$$
b_{1}=\begin{pmatrix} b_{++} & 0 \\
                    0 & b_{--}
    \end{pmatrix}
$$
in the Lie algebra
$\ug_{0}:=\ug_{1}\cap(\ug(\Hc_{+})\times\ug(\Hc_{-}))$ and an
element
$$
b_{2}=\begin{pmatrix} 0 & b_{+-} \\
                     b_{-+} & 0
    \end{pmatrix}$$
in the topological complement $\mathfrak{m}=
\ug(\Hc)\cap\left(\Sg_{2}(\Hc_{+},
\Hc_{-})\oplus\Sg_{2}(\Hc_{-},\Hc_{+})\right)$ of $\ug_{0}$ in
$(\ug_{\res})_{*}$. Accordingly,
$$\begin{array}{ll}
\|\alpha_{t}(b)\|_{(\ug_{\res})_{*}}& =
\|\exp(t\rho)b\exp(-t\rho)\|_{(\ug_{\res})_{*}}\\& =
\|\exp(t\rho)b_{1}\exp(-t\rho) +
\exp(t\rho)b_{2}\exp(-t\rho)\|_{(\ug_{\res})_{*}} \\ & =
\|e^{\ad(t\rho)}(b_{1}) +
e^{\ad(t\rho)}(b_{2})\|_{(\ug_{\res})_{*}}.
\end{array}
$$
Since $\ad(t\rho)$ preserves both $\ug_{0}$ and $\mathfrak{m}$, it
follows that
$$e^{\ad(t\rho)}(b_{1}) \in \ug_{0} ~~~~\text{ and }~~~~
e^{\ad(t\rho)}(b_{2})\in \mathfrak{m}.$$ By the very definition of
the norm $\|\cdot\|_{(\ug_{\res})_{*}}$, one  has
$$
\|\alpha_{t}(b)\|_{(\ug_{\res})_{*}} =
\|e^{\ad(t\rho)}(b_{1})\|_{1} + \|e^{\ad(t\rho)}(b_{2})\|_{2},
$$
where $\|\cdot\|_{1}$ (respectively $\|\cdot\|_{2}$) is the usual
norm in $\Sg_{1}$ (respectively $\Sg_{2}$). Since the conjugation
by a unitary element preserves both $\|\cdot\|_{1}$ and
$\|\cdot\|_{2}$, it follows that $\alpha_{t}$ acts by isometries
on $(\ug_{\res})_{*}$, in particular $\sup\limits_{t\in{\mathbb
R}}\Vert\alpha_t\Vert<\infty$.
\end{proof}

\begin{remark}\label{isometry}
\normalfont
The calculation in the proof of Proposition~\ref{block_diagonal}
actually shows that for every $u\in\U_{\res}$ satisfying $[d,u]=0$
we have $\Vert ubu^{-1}\Vert_{\res}=\Vert b\Vert_{\res}$
whenever $b\in\Bc_{\res}$.
In fact
$$\Vert ubu^{-1}\Vert_{\res}
=\Vert ubu^{-1}\Vert+\Vert[d,ubu^{-1}]\Vert_2 =\Vert b\Vert+\Vert
u[d,b]u^{-1}\Vert_2 =\Vert b\Vert+\Vert [d,b]\Vert_2 =\Vert
b\Vert_{\res}$$ where the second equality follows since $[d,u]=0$.
\end{remark}

\begin{corollary}\label{finite_rank}
If $\rho\in(\ug_{\res})_*$ is a finite-rank operator,
then the coadjoint isotropy group of $\rho$ is a Banach Lie subgroup
of $\U_{\res}$ and the corresponding $\U_{\res}$-coadjoint orbit $\Oc_\rho$
is a smooth leaf of the characteristic distribution of $(\ug_{\res})_*$.
\end{corollary}

\begin{proof}
The set of finite-rank operators $\Fc$ is a
dense subset of the predual $(\ug_{\res})_{*}$.
For every skew-symmetric finite-rank
operator $F$ there exists a unitary operator $u\in\1+\Fc$,
 such that $uFu^{-1}$ leaves both $\Hc_{-}$ and $\Hc_{+}$
 invariant.
(This follows since any two finite-rank operators are contained in
a certain finite-dimensional Lie algebra of finite-rank operators;
see for instance Lemma~1 in Chapter~I of \cite{dlH72} or
Proposition~3.1 in \cite{St75}.)
Note that $u\in \U_{\res}$,
and the isotropy groups of the elements $F$ and
$uFu^{-1}$ are conjugated by the element $u$.
Hence the isotropy
group at any finite-rank operator is a Banach-Lie subgroup of
$\U_{\res}$, and this shows that the conclusion of
Proposition~\ref{block_diagonal} is satisfied if we replace the hypothesis
$[d,\rho]=0$ by the condition that $\rho$ is a finite-rank operator.
\end{proof}

\begin{corollary}\label{hinkkanen}
Assume that $\rho\in(\ug_{\res})_*$  and that
there exist an orthonormal basis $\{e_n\}_{n\ge1}$
and the real numbers $t\in(0,1)$ and $s\in(0,3(1-t)/100]$
such that the following conditions are satisfied:
\begin{itemize}
\item[{\rm(i)}]
  We have $\{e_n\mid n\ge1\}\subseteq\Hc_{+}\cup\Hc_{-}$.
\item[{\rm(ii)}]
  The matrix $(\rho_{mn})_{m,n\ge1}$ of $\rho$ with respect to the basis
   $\{e_n\}_{n\ge1}$ has the properties
$$\vert\rho_{m+1,n+1}\vert\le t\vert\rho_{m,n}\vert\text{ whenever }m,n\ge1,$$
and
$$\vert\rho_{m,n}\vert^2\le\frac{s^2}{(mn)^2}\vert \rho_{mm}\rho_{nn}\vert
\text{ whenever }m,n\ge1\text{ and }m\ne n.$$
\end{itemize}
Then the coadjoint isotropy group of $\rho$ is a Banach-Lie subgroup
of $\U_{\res}$ and the corresponding $\U_{\res}$-coadjoint orbit $\Oc_\rho$
is a smooth leaf of the characteristic distribution of $(\ug_{\res})_*$.
\end{corollary}

\begin{proof}
It follows at once by Theorem~1 in \cite{Hk85} that there exists
an operator $a=-a^*\in\Sg_2(\Hc)$ such that the operator $u \rho
u^{-1}$ is diagonal with respect to the basis $\{e_n\}_{n\ge1}$,
where $u=\exp a$. In particular we have
$u\in\U_2\subseteq\U_{\res}$ and $[d,u \rho u^{-1}]=0$, so that we
can use Proposition~\ref{block_diagonal} to get the desired
conclusion.
\end{proof}

\begin{remark}\label{carey}
\normalfont Let $\rho\in\Bc(\Hc)$. In addition to the applications
of Proposition~\ref{block_diagonal} in the proofs of Corollaries
\ref{finite_rank}~and~\ref{hinkkanen}, we note that each of the
following two conditions is equivalent to the existence of an
unitary operator $u\in\U_{\res}$ such that $[d,u \rho u^{-1}]=0$:
\begin{itemize}
\item[{\rm(i)}]
There exists $p\in\Bc(\Hc)$ such that $p=p^*=p^2$, $p-p_{+}\in\Sg_2(\Hc)$,
and $\rho p=p\rho$.
\item[{\rm(ii)}]
There exists an element $\Wc\in\Gr$ such that $\rho(\Wc)\subseteq\Wc$.
\end{itemize}
In fact, our assertion concerning~(i) follows at once since
$$\{p\in\Bc(\Hc)\mid p=p^*=p^2\text{ and }p-p_{+}\in\Sg_2(\Hc)\}
=\{up_{+}u^{-1}\mid u\in\U_{\res}\}$$
according to Lemma~3.1 in \cite{Ca85}.

On the other hand, the assertion on condition~(ii) holds since
by Proposition~7.1.3 in \cite{PS90} we have
$$\Gr=\{u(\Hc_{+})\mid u\in\U_{\res}\}$$
and, in addition, if $p\in\Bc(\Hc)$ is the orthogonal projection
onto some closed subspace $\Wc\subseteq\Hc$ then
$\rho(\Wc)\subseteq\Wc$ if and only if $[p,\rho]=0$.
\end{remark}

\section{Some smooth adjoint orbits of
the restricted unitary  group}\label{adj}

For the sake of completeness, we are going to investigate in this section
the smoothness of adjoint orbits of the restricted unitary group.
In particular, we shall find sufficiently many smooth adjoint orbits 
of $\U_{\res}$
to fill an open subset of the Lie algebra $\ug_{\res}$
(Proposition~\ref{adjoint} below).

\begin{lemma}\label{riccati}
Assume that the element
$$\rho=\begin{pmatrix} \rho_{++} & \rho_{+-} \\
                     \rho_{-+} & \rho_{--}
    \end{pmatrix}\in\ug_{\res}$$
satisfies the conditions
\begin{equation}\label{disjoint}
\sigma(\rho_{++})\cap\sigma(\rho_{--})=\varnothing,
\end{equation}
and
\begin{equation}\label{small}
\Vert\rho_{+-}\Vert_2<\frac{1}{2}\dist(\sigma(\rho_{++}),\sigma(\rho_{--})).
\end{equation}
Then there exists $u\in\U_{\res}$
such that $[d,u^{-1}\rho u]=0$.
\end{lemma}

\begin{proof}
The hypotheses \eqref{disjoint}~and~\eqref{small}
imply that there exists a Hilbert-Schmidt operator
$k\colon\Hc_{+}\to\Hc_{-}$ satisfying the operator Riccati equation
$$k\rho_{+-}k+k\rho_{++}-\rho_{--}k=\rho_{-+}.$$
(This result was obtained in \cite{Mo95};
see also Theorem~4.6 and Remark~4.7 in \cite{ALT01},
as well as \cite{AMM03}.)
Then the operator
$$g=\begin{pmatrix}
    \id_{\Hc_{+}} & k^*           \\
                k & -\id_{\Hc_{-}}
    \end{pmatrix} $$
is invertible and has the properties
$[d,g]\in\Sg_2(\Hc)$, $g=g^*$, $[d,g^2]=0$ and
\begin{equation}\label{almost}
[d,g^{-1}\rho g]=0
\end{equation}
(see Subsection~2.3 in \cite{ALT01}).
Now let $g=us$ be the polar decomposition of the invertible operator
$g\in\Bc(\Hc)$, where $u\in\Bc(\Hc)$ is unitary and $s=(g^*g)^{1/2}$.

On the other hand, since $d^*=-d$,
it follows that the commutant $\{d\}'$ is a von Neumann
algebra of operators on $\Hc$.
Thus, since $g=g^*$ and $g^*g=g^2\in\{d\}'$, it is straightforward
to deduce that $(g^*g)^{1/2}\in\{d\}'$, that is, $[d,s]=0$.
Now recall that $[d,g]\in\Sg_2(\Hc)$ to deduce that
the unitary operator $u=gs^{-1}$ satisfies $[d,u]\in\Sg_2(\Hc)$,
that is, $u\in\U_{\res}$.

Moreover by \eqref{almost} we have
$$0=[d,g^{-1}\rho g]=[d,s^{-1}u^{-1}\rho us]=s^{-1}[d,u^{-1}\rho u]s,$$
where the latter equality follows since we have seen that $[d,s]=0$.
Now we get $[d,u^{-1}\rho u]=0$, as desired.
\end{proof}

\begin{proposition}\label{adjoint}
There exists an open $\U_{\res}$-invariant  neighborhood $V$ of
$d\in\ug_{\res}$ such that $V$ is a union of smooth adjoint orbits
of the Banach Lie group~$\U_{\res}$.
\end{proposition}

\begin{proof}
Denote by $V_0$ the set of all elements
$$\rho=\begin{pmatrix} \rho_{++} & \rho_{+-} \\
                     \rho_{-+} & \rho_{--}
    \end{pmatrix}\in\ug_{\res}$$
satisfying conditions
$$\sigma(\rho_{\pm\pm})\subseteq
\{y\in\ie{\mathbb R}\mid \vert y\mp\ie\vert<1/3\}$$
and
$$\Vert\rho_{\pm}\Vert_2<\frac{2}{3}.$$
It is clear that $V_0$ is an open neighborhood of $d\in\ug_{\res}$.
We are going to show that the set
$$V:=\bigcup_{u\in\U_{\res}}\Ad_{\U_{\res}}(u)V_0\subseteq\ug_{\res}$$
has the desired properties.

Indeed, $V$ is clearly invariant under the adjoint action of $\U_{\res}$,
it is a union of open sets, and one of these open sets contains $d$.
Moreover, it follows by Lemma~\ref{riccati} along with
the construction of $V$ that for every $\rho\in V$
there exists $u\in\U_{\res}$ such that $[d,u^{-1}\rho u]=0$.
Next denote $\tilde{\rho}=u^{-1}\rho u$, so that
$\exp(t\rho)=u\exp(t\tilde{\rho})u^{-1}$ for all $t\in{\mathbb R}$.
Then for all $t\in{\mathbb R}$ and $b\in\ug_{\res}$
it follows by means of Remark~\ref{isometry} that
$$\begin{aligned}
\Vert \exp(t\rho)\cdot b\cdot\exp(-t\rho)\Vert_{\res}
& =\Vert u\exp(t\tilde{\rho})\cdot u^{-1}\cdot b\cdot u\cdot
\exp(-t\tilde{\rho})\cdot u^{-1}\Vert_{\res} \\
& \le\Vert u\Vert_{\res}\cdot
\Vert \exp(t\tilde{\rho})\cdot u^{-1}\cdot b\cdot u\cdot
\exp(-t\tilde{\rho})\Vert_{\res}
\cdot \Vert u^{-1}\Vert_{\res} \\
&=\Vert u\Vert_{\res}\cdot
\Vert u^{-1}bu\Vert_{\res}
\cdot \Vert u^{-1}\Vert_{\res} \\
& \le\Vert u\Vert_{\res}^2\cdot\Vert u^{-1}\Vert^2_{\res}\cdot
\Vert b\Vert_{\res}.
\end{aligned}$$
Consequently the 1-parameter group
$$\alpha\colon{\mathbb R}\to\Bc(\ug_{\res}),\quad
\alpha_t(b)=\exp(t\rho)\cdot b\cdot \exp(-t\rho)$$
satisfies
$$\sup\limits_{t\in{\mathbb R}}\Vert\alpha_t\Vert
\le\Vert u\Vert_{\res}^2\cdot\Vert u^{-1}\Vert^2_{\res}.$$
Now the arguments in Remark~\ref{attempt} show that
the adjoint isotropy group of $\rho$ is a Lie subgroup of $\U_{\res}$,
and thus the adjoint orbit of $\rho$ is smooth.
\end{proof}

\begin{corollary}\label{coadjoint}
There exists an open $\U_{1,2}$-invariant open neighborhood
$V$ of $d\in\ug_{\res}=\ug_{1,2}^*$ such that
$V$ is a union of smooth coadjoint orbits
of the Banach-Lie group~$\U_{1,2}$.
\end{corollary}

\begin{proof}
Apply Proposition~\ref{adjoint} along with the fact that
$\U_{1,2}\hookrightarrow\U_{\res}$ and
the adjoint action of $\U_{\res}$ restricts to the coadjoint action
of~$\U_{1,2}$.
\end{proof}

\section{The Banach Lie-Poisson space associated to the central
extension of  $\ug_{2}$}\label{MAIN}

Denote by $\tilde{\ug}_{2} :=\ug_{2} \oplus \R$ the central
extension of $\ug_{2}$ defined by the restriction of $s$ to
$\ug_{2} \times \ug_{2}$,
where $s$ is the two-cocycle defined in \eqref{schwinger}.
The natural isomorphism
$(\tilde{\ug}_{2})^* \simeq \tilde{\ug}_{2}$ implies that
$\tilde{\ug}_{2}$ is a Banach Lie-Poisson space, for the Poisson
bracket given by
$$
\{f,g\}_{d}(\mu, \gamma) := \langle \mu,  \left[ D_{\mu}f(\mu),
D_{\mu}g(\mu) \right] \rangle - \gamma s(D_{\mu}f, D_{\mu}g)
$$
where $f, g \in \Ci(\tilde{\ug}_{2})$, $(\mu, \gamma)$ is
an arbitrary element in $\tilde{\ug}_{2}$, and $D_{\mu}$ denotes
the partial Fr\'echet derivative with respect to $\mu \in
\ug_{2}$.

\begin{theorem}\label{gr_leaf_in_u2}
The characteristic distribution of the Banach Lie-Poisson space
$\tilde{\ug}_{2}$ is integrable.
\end{theorem}

\begin{proof}
In order to prove that the characteristic distribution is
integrable, it suffices to check that all of the affine coadjoint
isotropy groups are Lie subgroups of the Hilbert Lie group~$\U_2$.
For this purpose we note that, for arbitrary
$(\mu,\gamma)\in\tilde{\ug}_2$, the corresponding isotropy group
of the affine coadjoint action of $\U_2$ on $\tilde{\ug}_2$ is
$$(\U_2)_{(\mu,\gamma)}
=\{g\in\U_2\mid \mu=g\mu g^{-1}-\gamma gdg^{-1}+\gamma d\},$$
according to the explicit expression of the affine coadjoint
action in Proposition~\ref{affine_coadjoint}. The previous
equality implies that
$$(\U_2)_{(\mu,\gamma)}
=\{g\in{\mathbb C}\1+\Sg_2(\Hc)\mid g^*g=gg^*=1\text{ and }
\mu=g\mu g^{-1}-\gamma gdg^{-1}+\gamma d\},$$ and now it is clear
that $(\U_2)_{(\mu,\gamma)}$ is an algebraic subgroup of
degree~$\le2$ of the group of invertible elements in the unital
Banach algebra ${\mathbb C}\1+\Sg_2(\Hc)$. Then the Harris-Kaup
theorem (see for instance Theorem~4.13 in \cite{Be05}) implies
that $(\U_2)_{(\mu,\gamma)}$ is a Lie group with respect to the
topology inherited from ${\mathbb C}\1+\Sg_2(\Hc)$. In particular,
this topology coincides with the one inherited from  $\U_2$. Since
$\U_2$ is a Hilbert Lie group, hence the Lie algebra of
$(\U_2)_{(\mu,\gamma)}$ has a complement in the Lie algebra of
$\U_2$, it then follows that $(\U_2)_{(\mu,\gamma)}$ is a
Banach Lie subgroup of $\U_2$, and this concludes the proof.
(Compare Remark~\ref{attempt}.)
\end{proof}

The transitivity of the action of the Lie group $\U_2$ on the
 connected component $\Gr^0$ of the
restricted Grassmannian has been established in Theorem 3.5 in
\cite{Ca85}, and Proposition V.7 in \cite{Ne02a}. 
That the
action of the subgroup $\U_{1,2}$ of $\U_2$ on $\Gr^0$ is transitive 
has been
proved in section~1.3.4 of \cite{Tu05} with the help of the canonical
basis defined in section~7.3 of \cite{PS90} and associated to any
element of the restricted Grassmannian. 
Below we give a shorter
and geometrical proof of the latter fact.

\begin{proposition}\label{homogeneous_under_u2}
The connected component $\Gr^0$ of the restricted Grassmannian is
a homogeneous space under the unitary group
$\U_{1,2}\subset\U_{2}$.
\end{proposition}

\begin{proof}
The restricted Grassmannian is a  symmetric space of the
restricted unitary group $\U_{\res}$. It follows from the
description of geodesics in Proposition 8.8 in \cite{Ar03} (see
also \cite{ON83} and \cite{CE75} or its infinite-dimensional
version as given in Example 3.9 in \cite{Ne02c}, or
Proposition~1.9 in \cite{Tu06}) that each geodesic
 of $\Gr$ starting at $W\in\Gr^0$ is given by
\begin{equation}\label{geodesic}
\beta(t)= (\exp tX)\cdot\Hc_{+},\quad X\in\mathfrak{m}_{W},
\end{equation}
where $\mathfrak{m}_{W}$ is the orthogonal in $\ug_{\res}$ to the
Lie algebra of the isotropy group of $W$. For $W=\Hc_+$ we have
$\mathfrak{m}=\ug(\Hc)\cap\left(\Sg_{2}(\Hc_{+},
\Hc_{-})\oplus\Sg_{2}(\Hc_{-},\Hc_{+})\right)$, and for
$W=g\cdot\Hc_{+}$ with $g\in\U_{\res}$, we have
$\mathfrak{m}_{W}=g\,\mathfrak{m}\,g^{-1}$. Note that for
$X\in\mathfrak{m}$, $\exp tX$ belongs to $\U_{1,2}\subset\U_{2}$.
Since the Hopf-Rinow Theorem is no longer true in the infinite
dimensional case, it is not clear whether every two elements in
the complete connected manifold $\Gr^{0}$ can be joined by a
geodesic. Nevertheless Theorem~B in \cite{Ek78} asserts that, for
every $W\in\Gr^0$, the set of elements which can be joined to $W$
by a unique minimal geodesic contains a dense $G_{\delta}$ set.
Moreover from the properties of the Riemannian exponential map,
there exists a neighborhood $\mathcal{V}$ of $\Hc_{+}$ in
$\Gr^{0}$ such that every element in $\mathcal{V}$ can be joined
to $\Hc_{+}$ be a (minimal) geodesic. Hence an arbitrary element
$W\in\Gr^0$ can be joined to an element $W'\in\mathcal{V}$ by a
geodesic
$$
\beta_1(t)= (\exp tX_1)\cdot W',\quad
X_1\in\mathfrak{m}_{W'},\quad t\in[0\,,1],
$$
and $W'$ can be joined to $\Hc_+$ by a geodesic
$$
\beta_2(t)= (\exp tX_2)\cdot\Hc_{+},\quad X_2\in\mathfrak{m},\quad
t\in[0\,,1].
$$
Consequently
$$
W=\beta_1(1)=(\exp X_1)\cdot W'=(\exp X_1)(\exp X_2)\cdot\Hc_+.
$$
But $X_1$ belongs to
$\mathfrak{m}_{W'}=\exp(X_2)\mathfrak{m}\exp(-X_2)$, hence
$$
W =(\exp X_2 \exp X_3)\cdot\Hc_+
$$
where $X_3=\textrm{Ad}\left(\exp(-X_2)\right)(X_1)$ belongs to
$\mathfrak{m}$. Since $\exp X_3$ and $\exp X_2$ are elements of
the unitary group $\U_{1,2}$, it follows that their product
belongs to $\U_{1,2}$. Thus $\U_{1,2}$ acts transitively on
$\Gr^0$.
\end{proof}

\begin{theorem}\label{main}
The connected component $\Gr^0$ of the restricted Grassmannian is
a strong symplectic leaf in the Banach Lie-Poisson space
$\tilde{\ug}_{2}$. More precisely, for every $\gamma\neq0$, the
$\U_{2}$-affine coadjoint orbit $\tilde{\mathcal{O}}_{(0,\gamma)}$
of $(0, \gamma)\in\tilde{\ug}_{2}$ is diffeomorphic to $\Gr^{0}$
via the application
$$
\begin{aligned}
\Phi_{\gamma}\Gr^0& \to \mathcal{O}_{(0,\gamma)}\\
                W & \mapsto 2\ie\gamma(p_{W}-p_{+}),
\end{aligned}
$$
where $p_{W}$ denotes the orthogonal projection on $W$. The
pull-back by $\Phi_{\gamma}$ of the symplectic form on
$\tilde{\mathcal{O}}_{(0,\gamma)}$ is $(\gamma/2)$-times
the symplectic form $\omega_{\Grr}$ on $\Gr^0$.
\end{theorem}

\begin{proof}
The assertion follows by the method of proof of
Theorem~\ref{gr_leaf_in_predual},
since $\Gr^0$ is transitively acted upon by the group $\U_2$
according to Proposition~\ref{homogeneous_under_u2}.
\end{proof}

Next we shall investigate the existence of
invariant complex structures
on certain covering spaces of the symplectic leaves of $\tilde{\ug}_2$
(Corollary~\ref{C3} below).
To this end we need two facts holding in a more general setting.
In connection with the first of these statements,
we note that invariant complex structures on
certain homogeneous spaces related to derivations of $L^*$-algebras
have been previously obtained by a different method in Theorem~IV.5
in \cite{Ne00}.

\begin{proposition}\label{P1}
Let $\Xg$ be a real Hilbert Lie algebra with a scalar product
denoted by $(\cdot\mid\cdot)$.
Assume that there exists a connected Hilbert Lie group $\U_{\Xg}$
whose Lie algebra is $\mathfrak{X}$; we write $\Lie(\U_{\Xg})=\Xg$.

Now let $D\colon\Xg\to\Xg$ be a bounded linear derivation
such that
\begin{equation}\label{skew}
(\forall x,y\in\Xg)\qquad (Dx\mid y)=-(x\mid Dy).
\end{equation}
Consider the closed subalgebra $\hg_0:=\Ker D$ of $\Xg$ and
define
$$H_0:=\langle\exp_{\U_{\Xg}}(\hg_0)\rangle,$$
that is, the subgroup of $\U_{\Xg}$ generated by the image of
$\hg_0$ by the exponential map.

If it happens that $H_0$ is a Lie subgroup of $\U_{\Xg}$, then
the smooth homogeneous space $\U_{\Xg}/H_0$ has
an invariant complex structure.
\end{proposition}

\begin{proof}
Denote $\Lg:=\Xg_{\mathbb C}$, that is,
the complex Hilbert-Lie algebra
which is the complexification of $\Xg$ and
is endowed with the complex scalar product $(\cdot\mid\cdot)$
extending the scalar product of $\Xg$.
We denote the complex linear extension of $D$ to $\Lg$ again by $D$.

Then $D^*=-D$ as operators on the complex Hilbert space $\Lg$,
so that $-\ie D\in\Bc(\Lg)$ is a self-adjoint operator.
Let us denote its spectral measure by $\delta\mapsto E(\delta)$.
Thus $E(\cdot)$ is a spectral measure on ${\mathbb R}$ and we have
$$D=\ie\int\limits_{\mathbb R}t\de E(t).$$
Also denote $S=(-\infty,0]$,
which is a closed subsemigroup of ${\mathbb R}$, and
$$\kg:=\Ran E(-S)=\Ran E([0,\infty))\subseteq\Lg.$$
Then $\kg$ is a closed subspace of $\Lg$
since it is the range of an idempotent continuous map.
In addition,
since $D$ is a derivation of the Hilbert Lie algebra $\Xg$
and $S$ is a closed semigroup,
it follows by Proposition~6.4 in \cite{Be05}
that $\kg$ is a complex subalgebra of $\Lg$
with the following properties:
\begin{itemize}
\item[{\rm(i)}] $[\hg_0,\kg]\subseteq\kg$,
\item[{\rm(ii)}] $\kg\cap\overline{\kg}=\hg_0+\ie\hg_0$ ($=\Ker D$), and
\item[{\rm(iii)}] $\kg+\overline{\kg}=\Lg$.
\end{itemize}
Moreover, for every $y\in\hg_0$ and all $x\in\Xg$ we have
$$D[y,x]=[Dy,x]+[y,Dx]=[y,Dx]$$
since $Dy=0$.
Therefore, we have $D\circ\ad_{\Xg}y=\ad_{\Xg}y\circ D$ for each $y\in\hg_0$.
According to the definition of $H_0$,
it then follows that for arbitrary $h\in H_0$
we have $\Ad_{\U_{\Xg}}h\circ D=D\circ\Ad_{\U_{\Xg}}h$ on $\Xg$.
Then the latter equality holds throughout $\Lg$, and it then follows
that the operator $\Ad_{\U_{\Xg}}h\colon\Lg\to\Lg$
commutes with every value of the spectral measure $E(\cdot)$.
In particular we have
$\Ad_{\U_{\Xg}}(h)\circ E(-S)=E(-S)\circ\Ad_{\U_{\Xg}}(h)$,
whence
\begin{itemize}
\item[{\rm(i')}]
$(\forall h\in H_0)\qquad\Ad_{\U_{\Xg}}(h)\kg\subseteq\kg$.
\end{itemize}
Now Theorem~6.1 in \cite{Be05} shows that
the smooth homogeneous space $\U_{\Xg}/H$
has an invariant complex structure.
\end{proof}

\begin{proposition}\label{P2}
Let $\Hc$ be an infinite-dimensional complex Hilbert space
and let $a\in\Bc(\Hc)$ such that $a^*=-a$.
Denote by
$$D=\ad_{\ug_2}\,a\colon\ug_2\to\ug_2,\quad x\mapsto[a,x]$$
the derivation of the compact $L^*$-algebra~$\ug_2$ defined by $a$,
and denote
$$\hg_0:=\Ker D=\{x\in\ug_2\mid [a,x]=0\}.$$
Next denote
$$H:=\{u\in\U_2\mid uau^{-1}=a\}$$
and in addition define
$$H_0:=\langle\exp(\hg_0)\rangle.$$
That is, $H_0$ is the subgroup of $\U_2$ generated by the image of
$\hg_0$ by the exponential map.
Then the following assertions hold:
\begin{itemize}
\item[{\rm(j)}] Both $H$ and $H_0$ are Lie subgroups of $\U_2$.
\item[{\rm(jj)}] The subgroup $H_0$ is the connected component of $\1\in H$.
\item[{\rm(jjj)}] The natural map
$$\U_2/H_0\to\U_2/H,\quad uH_0\mapsto uH,$$
is an $\U_2$-equivariant smooth covering map.
\end{itemize}
\end{proposition}

\begin{proof}
Consider  the Banach algebra $\Ac:={\mathbb C}\1+\Sg_2(\Hc)$ and
denote by $\varphi\colon\Ac\to{\mathbb C}$ the continuous linear
functional uniquely defined by the conditions $\varphi(\1)=1$ and
$\Ker\varphi=\Sg_2(\Hc)$. Then we have
$$
H=\{u\in\Ac^\times\mid u^*u=uu^*=\1\text{ and }\varphi(u)=1\}
$$
hence $H$ is a Lie subgroup of $\Ac^\times$ by the Harris-Kaup theorem
(see for instance Theorem~4.13 in \cite{Be05}),
and in addition the Lie algebra of $H$ is
$$
\Lie(H)=\{x\in\Ac\mid x^*=-x\text{ and }xa=ax\}=\hg_0.
$$
On the other hand, $H_0$ has the structure of connected Lie group
such that the inclusion map $H_0\hookrightarrow\U_2$ is an immersion
and $\Lie(H_0)=\hg_0$.
(See for instance Theorem~3.5 in \cite{Be05} and its proof.)
Since $H_0\subseteq H$ and $\Lie(H_0)=\Lie(H)=\hg_0$,
it then follows that $H_0$ is the connected component of $\1\in H$.
This can be seen directly by Lie theoretic methods;
specifically, one just has to use the fact that
the exponential map of any Banach Lie group is a local diffeomorphism
at $0$.
An alternative approach is to
use the proof of Lie's second theorem by means of
the Frobenius theorem
(see for instance Theorem~5.4 in Chapter~VI of \cite{La01}).
According to that proof,
the connected group $H_0$ is the integral manifold through $\1$
corresponding to a smooth left-invariant integrable distribution
on $\U_2$ whose fiber at $\1$ is (the complemented closed Lie subalgebra)
$\hg_0$.
Now recall the universality property of
the integral leaves of integrable distributions
according to Theorem~4.2 in Chapter~VI of \cite{La01}
or, more generally, Theorem~4(iii) in \cite{Nu92},
which implies that the inclusion map $H_0\hookrightarrow H$ is smooth.
Then the wished-for property that $H_0$ is open in $H$ follows since
$H_0$ and $H$ have the same tangent space at $\1\in H_0\subseteq H$.

By either of these methods it follows that $H_0$
is an open subgroup of the Lie subgroup $H$ of $\U_2$,
and then $H_0$ is in turn a Lie subgroup of $\U_2$.
Thus assertions (j) and (jj) are proved.
Assertion~(jjj) follows since the natural map
$\U_2/H_0\to\U_2/H$
is clearly an $\U_2$-equivariant map whose tangent
map at every point is an isomorphism.
\end{proof}

\begin{corollary}\label{C3}
Every symplectic leaf of the Hilbert Lie-Poisson space
$\tilde{\ug}_2$ is transitively acted on by $\U_2$ by means of the
affine coadjoint action and is $\U_2$-equivariantly covered by
some complex homogeneous space of $\U_2$.
\end{corollary}

\begin{proof}
Let $(\mu,\gamma)\in\tilde{\ug}_2$ arbitrary and denote
$a:=\mu-\gamma d\in\Bc(\Hc)$. With the notation of
Proposition~\ref{P2}, it is clear that $H$ is equal to the
isotropy group of the affine coadjoint action of $\U_2$. Thus the
symplectic leaf $\tilde{\Oc}_{(\mu,\gamma)}$ through
$(\mu,\gamma)$ is $\U_2$-equivariantly diffeomorphic to $\U_2/H$.
Now the conclusion follows since $\U_2/H$ is $\U_2$-equivariantly
covered by the complex homogeneous space $\U_2/H_0$, according to
Propositions \ref{P1}~and~\ref{P2}.
\end{proof}

\begin{remark}\label{R4}
\normalfont It follows by Corollary~\ref{C3} that every simply
connected symplectic leaf of the Banach Lie-Poisson space
$\tilde{\ug}_2$ has an $\U_2$-invariant complex structure. For
instance, this is the case for the connected component $\Gr^{0}$
of the restricted Grassmannian viewed as a symplectic leaf of
$\tilde{\ug}_2$ by means of Theorem~\ref{main}.
\end{remark}


\section{Some pathological properties of the restricted algebras}\label{pat}

\subsection{Unbounded unitary groups in the restricted algebra}
We are going to point out a property that
provides a good illustration for the difference between
the Banach $*$-algebra $\Bc_{\res}$ and a $C^*$-algebra
(Proposition~\ref{unbounded} below).

\begin{lemma}\label{exp}
Let $a\in\Bc(\Hc_{-},\Hc_{+})$ and assume that $a=v\vert a\vert$
and $a^*=w\vert a^*\vert$ are the polar decompositions of $a$ and
$a^*$, where $\vert a\vert\in\Bc(\Hc_{-})$ and $\vert
a^*\vert\in\Bc(\Hc_{+})$, while $v\colon\Hc_{-}\to\Hc_{+}$ and
$w\colon\Hc_{+}\to\Hc_{-}$ are partial isometries. Next, denote
$$\rho=\begin{pmatrix}    0 & a \\
                       -a^* & 0
       \end{pmatrix}\in\Bc(\Hc).
$$
Then
$$\exp\rho=\begin{pmatrix} \cos\vert a^*\vert & v\sin\vert a\vert \\
                          -w\sin\vert a^*\vert & \cos\vert a\vert
\end{pmatrix}.$$
\end{lemma}

\begin{proof}
We have
$$\rho^2=\begin{pmatrix} -aa^* &     0    \\
                          0    & -a^*a
         \end{pmatrix}
=-\begin{pmatrix} \vert a^*\vert^2 &          0    \\
                              0    & \vert a\vert^2
         \end{pmatrix}$$
hence
$$(\forall n\ge0) \quad \rho^{2n}=
(-1)^n\begin{pmatrix} \vert a^*\vert^{2n} &              0    \\
                                     0    & \vert a\vert^{2n}
         \end{pmatrix}.
$$
This implies that for every $n\ge0$ we have
$$\rho^{2n+1}=\rho\cdot\rho^{2n}
=(-1)^n \begin{pmatrix} 0 & v\vert a\vert \\
                       -w\vert a^*\vert & 0
\end{pmatrix}
\begin{pmatrix} \vert a^*\vert^{2n} &              0    \\
                                     0    & \vert a\vert^{2n}
         \end{pmatrix}
=(-1)^n \begin{pmatrix} 0 & v\vert a\vert^{2n+1} \\
                       -w\vert a^*\vert^{2n+1} & 0
\end{pmatrix}.
$$
Consequently
$$\exp\rho
=\sum_{n=0}^\infty\left(\frac{1}{(2n)!}\rho^{2n}
+\frac{1}{(2n+1)!}\rho^{2n+1}\right)
=\begin{pmatrix} \cos\vert a^*\vert & v\sin\vert a\vert \\
                          -w\sin\vert a^*\vert & \cos\vert a\vert
\end{pmatrix}$$
which concludes the proof.
\end{proof}

\begin{proposition}\label{unbounded}
All of the unitary groups $(\1+\Fc)\cap\U(\Hc)$, $\U_{1,2}$, and
$\U_{\res}$ are unbounded subsets of the unital associative Banach
algebra $\Bc_{\res}$.
\end{proposition}

\begin{proof}
We have
$$(\1+\Fc)\cap\U(\Hc)\subseteq\U_{1,2}\subseteq\U_{\res}$$
so it suffices to show that
\begin{equation}\label{unb}
\sup\{\Vert u\Vert_{\res}\mid u\in(\1+\Fc)\cap\U(\Hc)\}=\infty.
\end{equation}
To this end let $n\ge1$ be an arbitrary positive integer, pick a
projection $q_n=q_n^*=q_n^2\in\Bc(\Hc_{-})$ with $\dim(\Ran
q_n)=n$ and define
$a_n:=v_n((\pi/2)q_n)=(\pi/2)v_n\in\Bc(\Hc_{-},\Hc_{+})$, where
$v_n\colon\Hc_{-}\to\Hc_{+}$ is an arbitrary partial isometry such
that $v_n^*v_n=q_n$. Then $\vert a_n\vert=(\pi/2)q_n$, so that
$\sin\vert a_n\vert=q_n$ and then $\Vert(\sin\vert
a_n\vert)\Vert_2=\sqrt{\dim(\Ran q_n)}=\sqrt{n}$. Now
Lemma~\ref{exp} shows that the element
$$\rho_n=\begin{pmatrix}    0 & a_n \\
                       -a_n^* & 0
       \end{pmatrix}\in\ug(\Hc)\cap\Fc $$
satisfies
$$\Vert\exp(\rho_n)\Vert_{\res}\ge\Vert(\sin\vert a_n\vert)\Vert_2=\sqrt{n}.$$
Now the desired conclusion~\eqref{unb} follows since
$\exp(\rho_n)\in(\1+\Fc)\cap\U(\Hc)$ and $n\ge1$ is arbitrary.
\end{proof}

\subsection{The predual of the restricted algebra is not spanned by its
positive cone}
It is well known that every self-adjoint normal functional
in the predual of a $W^*$-algebra can be written as the difference
of two positive normal functionals.
It is also well known and easy to see that
a similar property holds for the preduals of numerous operator ideals.
More precisely,
if $\Jg$ and $\Bg$ are Banach operator ideals
such that the trace pairing
$$(\Bg,\Jg)\to{\mathbb C},\quad (T,S)\mapsto\Tr(TS)$$
is well defined and induces a topological isomorphism
of the topological dual $\Bg^*$ onto $\Jg$,
then for every $T=T^*\in\Bg$
there exist $T_1,T_2\in\Bg$ such that $T_1\ge0$, $T_2\ge 0$ and $T=T_1-T_2$.
In fact, we can take $T_1=(\vert T\vert+T)/2$ and $T_2=(\vert T\vert-T)/2$,
and we have $T_1,T_2\in\Bg$ since $\vert T\vert\in\Bg$.
(The latter property follows since if $T=W\vert T\vert$ is
the polar decomposition of $T$, then $\vert T\vert=W^*T\in\Bg$.)

We shall see in Proposition~\ref{order} below that the predual
$(\ug_{\res})_*$ of the restricted Lie algebra fails to have the
similar property of being spanned by its elements $\rho$ with
$\ie\rho\ge0$. In fact, the linear span of these elements turns
out to be the proper subspace $\ug_1$ of $(\ug_{\res})_*$.

\begin{lemma}\label{paulsen}
Let $\Hc_{\pm}$ be two complex separable Hilbert spaces,
$\Hc=\Hc_{+}\oplus\Hc_{-}$, $0\le a_{\pm}\in\Bc(\Hc_\pm)$, and
$t\in\Bc(\Hc_{-},\Hc_{+})$. Also denote
$$a=\begin{pmatrix}
a_{+} & t \\
t^* & a_{-}
\end{pmatrix}\in\Bc(\Hc).$$
Then the following assertions hold:
\begin{itemize}
\item[{\rm(i)}] We have $a\ge 0$ if and only if the inequality
\begin{equation}\label{lance}
\vert\langle\xi,t\eta\rangle\vert^2\le\langle\xi,a_{+}\xi\rangle\cdot
\langle\eta,a_{-}\eta\rangle
\end{equation}
holds for all $\xi\in\Hc_{+}$ and $\eta\in\Hc_{-}$.
\item[{\rm(ii)}] If $a\ge 0$ and in addition
$a_{\pm}\in\Sg_1(\Hc_{\pm})$ and $t\in\Sg_2(\Hc_{-},\Hc_{+})$,
then
\begin{equation}
\Vert t\Vert_2\le (\Tr a)/\sqrt{2}.
\end{equation}
\end{itemize}
\end{lemma}

\begin{proof}
For assertion~(i) see Exercise~3.2 at the end of Chapter~3 in
\cite{Pa02}.

Next, let $\{\xi_i\}_{i\ge1}$ and $\{\eta_j\}_{j\ge1}$ be
orthonormal bases in the Hilbert spaces $\Hc_{+}$ and $\Hc_{-}$,
respectively. Then \eqref{lance} shows that
$$(\forall i,j\ge1)\quad
\vert\langle\xi_i,t\eta_j\rangle\vert^2
\le\langle\xi_i,a_{+}\xi_i\rangle\cdot
\langle\eta_j,a_{-}\eta_j\rangle.$$ Now recall that $(\Vert
t\Vert_2)^2=\sum\limits_{i,j\ge1}\vert\langle\xi_i,t\eta_j\rangle\vert^2$,
$\Tr a_{+}=\sum\limits_{i\ge1}\langle\xi_i,a_{+}\xi_i\rangle$,
and $\Tr a_{-}=\sum\limits_{j\ge1}\langle\eta_j,a_{-}\eta_j\rangle$.
Thus, adding the above inequalities,
we get
$$(\Vert t\Vert_2)^2\le (\Tr a_{+})\cdot(\Tr a_{-})
\le (\Tr a_{+}+\Tr a_{-})^2/2 =(\Tr a)^2/2$$
and assertion~(ii) follows.
\end{proof}

\begin{proposition}\label{order}
The following assertions hold:
\begin{itemize}
\item[{\rm(i)}] If $a\in(\ug_{\res})_{*}$ and $\ie a\ge 0$, then
$a\in\Sg_1(\Hc)$ and $\Vert a\Vert_1\le\Vert
a\Vert_{(\ug_{\res})_{*}} \le(1+\sqrt{2})\Vert a\Vert_1$.
\item[{\rm(ii)}] If $\rho\in(\ug_{\res})_{*}\setminus\ug_1$ then
there exist no $\rho_1,\rho_2\in(\ug_{\res})_{*}$ such that
$\ie\rho_1\ge0$, $\ie\rho_2\ge0$, and $\rho=\rho_1-\rho_2$.
\end{itemize}
\end{proposition}

\begin{proof}
(i) Let $a\in(\ug_{\res})_{*}$ such that $\ie a\ge0$, and denote
$\ie a=:\begin{pmatrix}
a_{+} & t \\
t^* & a_{-}
\end{pmatrix}$.
Then
$$\begin{aligned}
\Vert a\Vert_1 &=\Vert\ie a\Vert_1=\Tr(\ie a) =\Tr a_{+}+\Tr a_{-}
=\Vert a_{+}\Vert_1+\Vert a_{-}\Vert_1 \\
&\le\Vert \ie a\Vert_{(\ug_{\res})_{*}} =\Vert a_{+}\Vert_1+\Vert
a_{-}\Vert_1
+2\Vert t\Vert_2 \\
&\le \Vert a_{+}\Vert_1+\Vert a_{-}\Vert_1 +\sqrt{2}\cdot\Tr (\ie a)
=(1+\sqrt{2})\Vert \ie a\Vert_1=(1+\sqrt{2})\Vert a\Vert_1 ,
\end{aligned}$$
where the second inequality follows by Lemma~\ref{lance}(ii).
Consequently, for all $a\in(\ug_{\res})_{*}$ with $\ie a\ge0$ we have
$\Vert a\Vert_1\le\Vert a\Vert_{(\ug_{\res})_{*}}
\le(1+\sqrt{2})\Vert a\Vert_1$.

(ii) Let $\rho\in(\ug_{\res})_{*}\setminus\ug_1$ and assume that
there exist elements $\rho_1,\rho_2\in(\ug_{\res})_{*}$ such that
$\ie\rho_1\ge0$, $\ie\rho_2\ge0$, and
$\rho=\rho_1-\rho_2$. Then $\ie\rho_1,\ie\rho_2\in\Sg_1(\Hc)$
according to the assertion~(i), which we have already proved.
Consequently, $\rho_1,\rho_2\in\ug_1$, whence
$\rho=\rho_1-\rho_2\in\ug_1$. This is a contradiction with the
assumption on $\rho$, which concludes the proof.
\end{proof}

\subsection{The Cartan subalgebras of $\ug_{\res}$ are
not $\U_{\res}$-conjugate}\label{cartan}

For a (finite-dimensional) compact connected semi-simple Lie
subgroup $G$ of the unitary group $\U(n)$, every element $X$ of
the Lie algebra $\mathfrak{g}$ of $G$ is conjugate to a diagonal
element by an element of $G$. This can be seen as follows (see
\cite{He62} chap.~V theorem 6.4 for more general results). Take a
diagonal element $H\in\mathfrak{g}$ such that the one-parameter
subgroup $\exp tH$ is dense in the torus whose Lie algebra is the
set of diagonal matrices belonging to $\g$. On $G$, consider the
continuous function $g\mapsto \textrm{B}\left(H,
\textrm{Ad}(g)(X)\right)$, where $\textrm{B}$ denotes the Killing
form of $G$. By compactness, this function takes a minimum at some
$g_{0}$, and for every element $Y$ in $\g$ one has
$$
\frac{d}{dt}{\textrm{B}\left(H,\Ad(\exp
tY)\Ad(g_{0})(X)\right)}|_{t=0} = 0,
$$
i.e $\textrm{B}\left(H, [Y, \Ad(g_{0})(X)]\right) = 0$. Since the
Killing form is $\Ad(G)$-invariant, one has
$$\textrm{B}\left(H, [Y, \Ad(g_{0})(X)]\right) =
\textrm{B}\left([\textrm{Ad}(g_{0})(X), H],Y\right).$$ The
non-degeneracy of the Killing form then implies that
$[\Ad(g_{0})(X), H] = 0$. But $H$ has been chosen such that the
centralizer of $H$ is the set of diagonal matrices belonging to
$\mathfrak{g}$. Consequently $\textrm{Ad}(g_{0})(X)$ is a diagonal
element in $\mathfrak{g}$.
It follows that the maximal Abelian
subalgebras, called {\it Cartan subalgebras}, of $\mathfrak{g}$ are
conjugate under $G$. Naturally this proof does not work anymore
for an infinite-dimensional group since the argument to minimize the
corresponding function is missing. In fact, we will show below that
the Cartan subalgebras of $\ug_{\res}$ are not
$\U_{\res}$-conjugate, in general.

We note that a related fact follows from results in the paper
\cite{BS96}. Specifically, let $\rho_0\in(\ug_{\res})_*$ such that
$[d,\rho_0]=0$, $\Ker\rho_0=\{0\}$, and each eigenvalue of
$\rho_0$ has multiplicity~1.
Next denote by $\Oc_{\rho_0}$ the
coadjoint $\U_{\res}$-orbit of $\rho_0$, let
$\rho\in(\ug_{\res})_*$, and define
$$f_{\rho}\colon\Oc_{\rho_0}\to(0,\infty),\quad
f_{\rho}(b)=\Vert\rho-b\Vert_2.$$ If the function $f_{\rho}$
happens to have a critical point $\rho_1\in\Oc_{\rho_0}$, then
$[\rho_1,\rho]=0$ according to \cite{BS96}. Since
$\rho_1\in\Oc_{\rho_0}$, there exists $u\in\U_{\res}$ such that
$\rho_1=u\rho_0u^{-1}$, and then $[\rho_0,u^{-1}\rho u]=0$. The
latter equality implies that $u^{-1}\rho u$ commutes with all of
the spectral projections of $\rho_0$. Hence $[d,u^{-1}\rho u]=0$
in view of the spectral assumptions on $\rho_0$, and then
Proposition~\ref{block_diagonal} applied to $u^{-1}\rho u$ shows
that the coadjoint isotropy group of $\rho$ is a Banach-Lie
subgroup of $\U_{\res}$ and the corresponding
$\U_{\res}$-coadjoint orbit $\Oc_\rho$ is a smooth leaf of the
characteristic distribution of $(\ug_{\res})_*$.

\begin{proposition}
The unitary group  $\U_{\res}$ does not act transitively on the
set of Cartan subalgebras of its Lie algebra.
\end{proposition}

\begin{proof}
Endow the Hilbert space $\Hc$ with an orthonormal basis $\mathcal{B} =
\{e_{n}\mid n\in \mathbb{Z}^{*}\}$, such that $\{e_{n}\mid n\in-\N^*\}$
is an orthonormal basis of $\Hc_+$ and $\{e_n\mid n\in\N^*\}$ an orthonormal
basis of $\Hc_-$. The set $\mathcal{D}$ of skew-Hermitian bounded
diagonal operators with respect to $\mathcal{B}$ form a Cartan
subalgebra of $\ug_{\res}$. Now consider the following  subset of
the set of anti-diagonal elements in $\ug_{\res}$~:
$$
\mathcal{J} = \{ J\in\ug_{\res}~|~J(e_{n})\in\R e_{-n}~~\forall n
\in \mathbb{Z}^*\}.
$$
Since the  coefficients $J_{-k,k}$, $k\in\mathbb{Z}^*$, of
$J\in\mathcal{J}$ satisfy $J_{-k,k} = -J_{k,-k}$, it follows from
an easy computation that $\mathcal{J}$ is Abelian. An element $B =
(B_{i,j})\in\ug_{\res}$ commutes with every element $J= (J_{i,j})$
in $\mathcal{J}$ if and only if
\begin{equation}\label{crochet}
([B, J]_{i,-k}) = (B_{i,k}J_{k,-k}-J_{i,-i}B_{-i,-k})
\end{equation}
vanishes for every $J \in \mathcal{J}$. This implies the following
conditions:
 $$\begin{aligned}
 B_{i,k}&= 0  ~\textrm{ for }~i\notin\{k, -k\};\\
 B_{k,k}&= B_{-k,-k}~\textrm{ for }k\in\mathbb{Z}^*;\\
B_{-k,k}&= -B_{k,-k}~\textrm{ for }k\in\mathbb{Z}^*.
\end{aligned}
$$
It follows that the maximal Abelian subalgebra $\mathcal{C}$ of
$\ug_{\res}$ which contains $\mathcal{J}$ is
$\mathcal{J}+\mathcal{D}_+$, where
$$\mathcal{D}_+=\{D=(D_{i,j})\in\mathcal{D}~|~D_{-k,-k} =
D_{k,k}~~\forall k\in\mathbb{Z}^*\}.
$$
Let us prove by contradiction that the Cartan subalgebras
$\mathcal{C}$ and $\mathcal{D}$ are not conjugate under
$\U_{\res}$. Suppose that there exists a unitary operator
$$g = \begin{pmatrix} g_{++} & g_{+-} \\
                     g_{-+} & g_{--}
    \end{pmatrix}
\in\U_{\res}$$ 
such that $g\mathcal{J}g^{-1}=\mathcal{D}$.
Consider an element $$J = \begin{pmatrix} 0 & J_{+-} \\
                     J_{-+} & 0\end{pmatrix}\in\mathcal{J}$$
which is a Hilbert-Schmidt
operator that is not trace class. 
One has
$$\begin{aligned}
gJg^{-1} &= \begin{pmatrix} g_{++} & g_{+-} \\
                            g_{-+} & g_{--}
            \end{pmatrix} 
            \begin{pmatrix} 0 & J_{+-} \\
                       J_{-+} & 0
            \end{pmatrix} 
            \begin{pmatrix} g_{++}^* & g_{-+}^* \\
                     g_{+-}^* & g_{--}^*
            \end{pmatrix}\\
&=\begin{pmatrix} g_{+-}J_{-+}g_{++}^* + g_{++}J_{+-}g_{+-}^* &
g_{+-}J_{-+}g_{-+}^* + g_{++}J_{+-}g_{--}^*\\
g_{--}J_{-+}g_{++}^*+g_{-+}J_{+-}g_{+-}^*&
g_{--}J_{-+}g_{-+}^*+g_{-+}J_{+-}g_{--}^*
\end{pmatrix}.
\end{aligned}
$$
By hypothesis, $gJg^{-1}$ is a diagonal operator 
$$D =  \begin{pmatrix} D_{++} & 0 \\
                     0 & D_{--}
    \end{pmatrix}
$$
with $D_{++} = g_{+-}J_{-+}g_{++}^* + g_{++}J_{+-}g_{+-}^*$ and
$D_{--} = g_{--}J_{-+}g_{-+}^*+g_{-+}J_{+-}g_{--}^*$.
    Now, since $g$
belongs to $\U_{\res}$, $g_{+-}$ and $g_{-+}$ are Hilbert-Schmidt.
Since $J$ belongs to $\Sg_2(\Hc)$, $J_{+-}$ and $J_{-+}$ are
Hilbert-Schmidt as well. From the relation
$\Sg_{2}\cdot\Sg_{2}\subset\Sg_1$, it follows that $D_{++}$ and
$D_{--}$ are trace class, hence $D$ belongs to $\Sg_{1}(\Hc)$. But
this implies that $J = g^{-1}Dg$ is also trace class, since
$\Sg_1(\Hc)$ is an ideal of $\mathcal{B}(\Hc)$. This leads to a
contradiction by the choice of $J\in\mathcal{J}$. It follows that
elements in $\mathcal{J}\setminus \Sg_{1}(\Hc)$ are not
$\U_{\res}$-conjugate to diagonal elements. Consequently, the
Cartan subalgebra $\mathcal{C}$ and $\mathcal{D}$ are not
$\U_{\res}$-conjugate.
\end{proof}

\begin{remark}\label{diffeos}
\normalfont
Since every skew-Hermitian operator is conjugate to a diagonal
operator by a unitary operator, the set of conjugacy classes of
Cartan subalgebras in $\ug_{\res}$ is in bijection with
$\U(\Hc)/\U_{\res}$ and is infinite. The conjugacy classes of
Cartan subalgebras are related to the conjugacy classes of maximal
tori. 
An infinite number of conjugacy classes of maximal tori has
already been encountered in the case of some groups of contactomorphisms 
(see \cite{Le01}). 
Examples of maximal tori of different dimensions
were provided in \cite{HT02} in some groups of
symplectomorphisms.
\end{remark}

\section*{Acknowledgments} 
The research of the first author was partially supported by Romanian
grant CEx05-D11-23/2005 and that of the second by the Swiss National Science Foundation. We thank Alan Weinstein who motivated this investigation by asking us whether the Grassmannians studied in 
\cite{Tu05} fall within the framework of the theory developed in  \cite{BR05}. Our thanks go also to Tilmann Wurzbacher for several discussions that influenced our presentation.


\begin{thebibliography}{AMM033}

\small

\bibitem[ALT01]{ALT01}
V.~Adamjan, H.~Langer, C.~Tretter,
Existence and uniqueness of contractive solutions of some Riccati equations,
\textit{J. Funct. Anal.} {\bf 179} (2001), no.~2, 448--473.

\bibitem[AMM03]{AMM03}
S.~Albeverio, K.A.~Makarov, A.K.~Motovilov,
Graph subspaces and the spectral shift function,
\textit{Canad. J. Math.} \textbf{55} (2003), no.~3, 449--503.

\bibitem[Ar03]{Ar03}
A.~Arvanitoyeorgos, {\it An Introduction to Lie Groups and the
Geometry of Homogeneous Spaces}, Student Math. Library no.~22,
American Math. Society, Providence, R.I., 2003.


\bibitem[Ba90]{Ba90}
A.G.~Baskakov,
Diagonalization of operators and complementability of
subspaces of Banach spaces,
\textit{Ukrain. Mat. Zh.} \textbf{42} (1990), no. 7, 867--873;
translation in
\textit{Ukrainian Math. J.} \textbf{42} (1990), no.~7, 763--768 (1991).

\bibitem[Be06]{Be05}
D.~Belti\c t\u a,
\textit{Smooth Homogeneous Structures in Operator Theory},
Chapman \& Hall/CRC Monographs and Surveys in Pure
   and Applied Mathematics, 137. Chapman \& Hall/CRC Press,
Boca Raton-London-New York-Singapore, 2006.

\bibitem[BP05]{BP05}
D.~Belti\c t\u a, B.~Prunaru,
Amenability, completely bounded projections,
dynamical systems and smooth orbits,
\textit{Integral Equations Operator Theory} 
(to appear). 
(See \textit{preprint} math.OA/0504313.)

\bibitem[BR05]{BR05}
D.~Belti\c t\u a, T.S.~Ratiu, Symplectic leaves in real Banach
Lie-Poisson spaces, \textit{Geom. Funct. Analysis} \textbf{15}
(2005), no.~4, 753--779.

\bibitem[BS96]{BS96}
R.~Bhatia, P.~\v Semrl,
Distance between Hermitian operators in Schatten classes,
\textit{Proc. Edinburgh Math. Soc. (2)} \textbf{39} (1996), no.~2, 377--380.

\bibitem[Ca85]{Ca85}
A.L.~Carey,
Some homogeneous spaces and representations of
the Hilbert Lie group ${\mathcal U}(H)\sb 2$,
\textit{Rev. Roumaine Math. Pures Appl.} {\bf 30} (1985), no.~7, 505--520.

\bibitem[CE75]{CE75}
J.~Cheeger, D.G.~Ebin, \textit{Comparison Theorems in Riemannian
Geometry}, North-Holland, Amsterdam, 1975.


\bibitem[Ek78]{Ek78}
I.~Ekeland, The Hopf-Rinow Theorem in infinite dimension,
\textit{J. Differential Geometry} {\bf 13} (1978), 287--301.


\bibitem[dlH72]{dlH72}
P.~de la Harpe,
\textit{Classical Banach-Lie Algebras and Banach-Lie Groups of Operators
in Hilbert Space}.
Lecture Notes in Mathematics, Vol.~285,
Springer-Verlag, Berlin-New York, 1972.

\bibitem[HT03]{HT02}
J.-C.~Hausmann, S.~Tolman, 
Maximal Hamiltonian tori for polygon spaces,  
\textit{Ann. Inst. Fourier (Grenoble)} 
\textbf{53} (2003), no.~6, 1925--1939.

\bibitem[He62]{He62} S.~Helgason,
\textit{Differential Geometry and Symmetric Spaces},
Academic Press, New York, 1962.

\bibitem[Hk85]{Hk85}
A.~Hinkkanen,
On the diagonalization of a certain class of operators,
\textit{Michigan Math. J.} \textbf{32} (1985), no.~3, 349--359.

\bibitem[La01]{La01}
S.~Lang,
\textit{Fundamentals of Differential Geometry}
(corrected second printing),
Graduate texts in mathematics, vol. 191,
Springer-Verlag,
New York,
 2001.


\bibitem[Le01]{Le01} 
E.~Lerman, 
Contact cuts, 
\textit{Israel J. Math.} \textbf{124} (2001), 77--92. 


\bibitem[Mo95]{Mo95}
A.K.~Motovilov,
Removal of the resolvent-like energy dependence from interactions and
invariant subspaces of a total Hamiltonian,
\textit{J. Math. Phys.} \textbf{36} (1995), no.~12, 6647--6664.


\bibitem[Ne00]{Ne00}
K.-H.~Neeb,
Highest weight representations and
    infinite-dimensional K\"ahler manifolds.
In: \textit{Recent Advances in Lie Theory (Vigo, 2000)},
Res. Exp. Math., 25, Heldermann,
    Lemgo, 2002, 367--392.


\bibitem[Ne02a]{Ne02a}
K.-H.~Neeb, Classical Hilbert-Lie groups, their extensions and
their homotopy groups. In: {\it Geometry and Analysis on Finite-
and Infinite-dimensional Lie Groups (B\c edlewo, 2000)}, Banach
Center Publ., 55, Polish Acad. Sci., Warsaw, 2002, 87--151.

\bibitem[Ne02b]{Ne02b}
K.-H.~Neeb, Universal central extensions of Lie groups. {\it Acta
Appl. Math.} {\bf 73} (2002), no.~1-2, 175--219.

\bibitem[Ne02c]{Ne02c}
K.-H.~Neeb, A Cartan-Hadamard theorem for Banach-Finsler
manifolds, {\it Geom. Dedicata} {\bf 95} (2002), 115--156.


\bibitem[Nu92]{Nu92}
F.~N\"ubel,
On integral manifolds for vector space distributions,
\textit{Math. Ann.} {\bf 294} (1992), no.~1, 1--17.


\bibitem[OR03]{OR03}
A.~Odzijewicz, T.S.~Ratiu,
Banach Lie-Poisson spaces and
reduction, \textit{Comm. Math. Phys.}
\textbf{243} (2003), no.~1, 1--54.

\bibitem[OR04]{OR04}
A.~Odzijewicz, T.S.~Ratiu,
Extensions of Banach Lie-Poisson spaces,
\textit{J. Funct. Anal.}
\textbf{217} (2004) no.~1, 103--125.

\bibitem[ON83]{ON83}
B.~O'Neill, \textit{Semi-Riemannian Geometry with Applications to
Relativity}, Academic Press, United Kingdom, 1983.

\bibitem[OrR04]{OrR04}
J.-P.~Ortega, T.S.~Ratiu,
\textit{Momentum Maps and Hamiltonian Reduction},
Progress in Mathematics, 222.
Birkh\"auser Boston, Inc., Boston, MA, 2004.



\bibitem[Pa02]{Pa02}
V.~Paulsen, {\it Completely Bounded Maps and Operator Algebras}.
Cambridge Studies in Advanced Mathematics, 78.
   Cambridge University Press, Cambridge, 2002.


\bibitem[PS90]{PS90}
A.~Pressley, G.~Segal,
\textit{Loop Groups},
Oxford Mathematical Monographs. Oxford Science Publications. The
   Clarendon Press, Oxford University Press, Oxford, 1990.

\bibitem[Se81]{Se81} G.~Segal, 
Unitary Representations of some infinite dimensional groups, 
\textit{Comm. Math. Phys} \textbf{80} (1981), no.~3, 301--342.

\bibitem[SW85]{SW85} G.~Segal, G.~Wilson, 
Loop groups and equations of KdV type, 
\textit{Publications Math\'ematiques de l'I.H.E.S.} \textbf{61} (1985),
5--6.


\bibitem[St75]{St75}
I.~Stewart,
\textit{Lie Algebras Generated by Finite-dimensional Ideals}.
Research Notes in Mathematics, Vol.~2,
Pitman Publishing, London-San Francisco, Calif.-Melbourne, 1975.


\bibitem[Tu05]{Tu05}
A.B.~Tumpach, \textit{Vari\'et\'es K\"ahl\'eriennes et
Hyperk\"ahl\'eriennes de Dimension Infinie}, Ph.D Thesis, 
\'Ecole Polytechnique, Paris, 2005.

\bibitem[Tu06]{Tu06} A.B.~Tumpach, 
Mostow Decomposition Theorem
for a $L^*$-group and applications to affine coadjoint orbits and
stable manifolds, \textit{preprint} math-ph/0605039 (May 2006).

\bibitem[Wu01]{Wu01}
T.~Wurzbacher,
Fermionic second quantization and the
    geometry of the restricted Grassmannian.
In \textit{Infinite Dimensional
    K\"ahler Manifolds (Oberwolfach, 1995)},
DMV Sem., 31, Birkh\"auser, Basel, 2001, 287--375.


\end{thebibliography}
\end{document}